\documentclass{article}
\usepackage{amssymb,amsmath,amsthm,latexsym,verbatim}
\newcounter{defcount}
\newcounter{thmcount}
\newtheorem{proposition}{Proposition}[section]
\newtheorem{lemma}[proposition]{Lemma}
\newtheorem{theorem}[thmcount]{Theorem}
\newtheorem*{thm}{Theorem}
\newtheorem{corollary}[thmcount]{Corollary}
\newtheorem{definition}[defcount]{Definition}
\theoremstyle{definition}
\newtheorem*{notation}{Notation}

\theoremstyle{definition}
\newtheorem{paragr}{}[section]

\begin{document}

\bibliographystyle{plain}

\title{Mixing on a Class of Rank One Transformations
        }

\author{Darren Creutz
\thanks{Darren.A.Creutz@williams.edu}
   \and C. E.  Silva
     \thanks{Dept. of Mathematics,
	Williams College, Williamstown, MA 01267.
	csilva@williams.edu} }

\date {August 10, 2001}

\maketitle

\begin{abstract}
We prove that a rank one transformation satisfying a condition called 
restricted growth is a mixing transformation if and only if the 
spacer sequence for the transformation is uniformly ergodic.  Uniform 
ergodicity is a generalization of the notion of ergodicity for 
sequences in the sense that the mean ergodic theorem holds for what we 
call dynamical sequences.  In particular, Adams' class of staircase 
transformations and Ornstein's class constructed using ``random 
spacers'' have both restricted growth and uniformly ergodic spacer sequences.
\end{abstract}

\medskip

\newpage

\thispagestyle{empty}
\clearpage
\setcounter{page}{1}

%%%%%%%%%%%%%%%%%%%%%%%%%%%%%%%%%%%%%%%%%%%%%%%%%%%%%%%%%%%%%%%%%%%%%%%%%%%%

\section{Introduction}

\begin{paragr}\textbf{Background.}
Rank one transformations are transformations defined inductively
using at each step a single Rohlin tower or column.  They gained 
importance with Ornstein's construction of a mixing transformation 
with no square root \cite{dO72}, and have been used as a source of 
examples and counterexamples in ergodic theory.  Rank one 
mixing transformations have been shown to enjoy many interesting 
properties.  As the culmination of several works it is now 
known that rank one mixing transformations are mixing of all orders 
and have the property of minimal self-joinings 
\cite{sK84}, \cite{jK88}, \cite{vR93}, and  in particular they are prime 
and have trivial centralizer.  We refer to \cite{nF70} for rank one 
constructions and \cite{sF97} for a recent survey of results.
\end{paragr}

\begin{paragr}\textbf{History.}
Ornstein's construction  of a rank 
one mixing transformation uses the notion of ``random spacers''
and yields a class of transformations 
so that almost surely a transformation in this class is
mixing.  Bourgain \cite{jB93}  
 showed that transformations in Ornstein's class 
almost surely have singular spectrum, and more recently El Houcein 
\cite{eA99} showed that transformations in Ornstein's class 
almost surely are mutually singular.  
However,  Ornstein's construction does not exhibit a specific rank one 
mixing transformation.  It was conjectured by Smorodinsky that a 
specific rank 
one transformation, now called the classical staircase transformation, 
is a 
mixing transformation.  Adams and Friedman \cite{AF92} (unpublished) constructed 
explicit mixing staircase  transformations, but the conjecture 
remained open until Adams 
\cite{tA98} showed that transformations in 
a class of staircase transformations, which 
includes the 
classical staircase, are mixing.  The classical staircase has 
been shown to have singular spectrum \cite{iK96}, and the construction and proof 
in \cite{tA98} have been generalized to mixing staircase $\mathbb 
Z^{d}$-actions in \cite{AS99}.
\end{paragr}

\begin{paragr}\textbf{Result.}
Our main result (Theorem \ref{T:main}) is that a rank one transformation, satisfying a 
condition we call restricted growth, is mixing if and only if 
its spacer sequence is uniformly ergodic.  A rank one transformation 
is specified by a sequence of positive integers called the cut 
sequence and a doubly-indexed sequence of integers, called the 
spacer sequence; this is a specific instance 
of what we call a dynamical sequence.  We generalize the notion of 
ergodicity of a sequence (as used in the Blum-Hanson theorem 
\cite{BH60}) to dynamical sequences, yielding notions of ergodicity 
and uniform ergodicity for dynamical sequences.  We introduce a 
condition called uniform mixing and show that it implies mixing.  The 
proof of our main theorem is then accomplished by showing that uniform 
ergodicity of the spacer sequence implies uniform mixing.
\end{paragr}

\begin{paragr}\textbf{Applications.}
We then apply our theorem to give another proof that the staircase 
transformations of Adams~\cite{tA98} are mixing.   For 
staircase transformations, restricted growth is equivalent to the 
condition under which Adams shows mixing.  Our technique for showing 
mixing is to first show that staircase transformations have mixing 
height sequences, then show that staircase transformations have 
uniformly ergodic spacer sequences using refinements of techniques in 
\cite{tA98}.  We conclude with another proof that Ornstein's class of 
transformations are almost surely mixing transformations \cite{dO72}.  This is 
accomplished using the notion of double ergodicity \cite{BFMS01} to 
show weak mixing, then using Ornstein's probabilistic lemma to show uniform 
ergodicity of the spacer sequence.
\end{paragr}

\section{Preliminaries}\label{S:prelims}

\begin{paragr}\textbf{Dynamical Systems.}
Let $(X,\mu)$ be a finite measure (probability) space isomorphic to the 
unit interval in $\mathbb{R}$ under Lebesgue measure.  We concern 
ourselves with invertible, measurable and measure-preserving 
transformations $T:X\to X$ on $(X,\mu)$ and use the term \textbf{transformation} to refer 
exclusively to such.  The space $(X,\mu)$ and a transformation $T$ on 
it form the \textbf{dynamical system} $(X,\mu,T)$.
\end{paragr}

\begin{paragr}\textbf{Ergodicity.}
A transformation $T: X\to X$ is 
\textbf{ergodic} when every T-invariant set is either null or 
conull---any measurable set $A$ such that $T(A) = A$ must have 
measure zero or full measure.
The \textbf{weak ergodic theorem} states that given an ergodic transformation 
$T$, for any measurable sets $A$ and $B$,
\[
\lim_{n\to\infty}\frac{1}{n}\sum_{i=0}^{n-1}\mu(T^{i}(A) \cap B) - 
\mu(A)\mu(B) = 0,
\]
while the \textbf{von Neumann 
(mean) ergodic theorem} guarantees 
convergence in the mean---for any measurable set $B$,
\[
\lim_{n\to\infty}\int\big{|}\frac{1}{n}\sum_{i=0}^{n-1}\chi_{B}\circ 
T^{-i} - \mu(B)\big{|}^{2} d\mu = 0.
\]
Given an ergodic transformation $T$, a 
sequence of integers $\{a_{n}\}$ is 
\textbf{ergodic with respect to $\boldsymbol{T}$} when for any 
measurable set $B$,
\[
\lim_{n\to\infty} \int\big{|}\frac{1}{n}\sum_{i=0}^{n-1}\chi_{B}\circ 
T^{-a_{i}} - \mu(B)\big{|}^{2} d\mu = 0
\]
and is \textbf{weak ergodic 
with respect to $\boldsymbol{T}$} when for any measurable sets $A$ 
and $B$, 
\[
\lim_{n\to\infty}\frac{1}{n}\sum_{i=0}^{n-1}\mu(T^{a_{i}}(A) \cap 
B) - \mu(A)\mu(B) = 0.
\]
Note that a sequence that is ergodic with respect to $T$ is weak 
ergodic with respect to $T$; the proof is left to the reader.  That 
the converse does not hold in general is due to Friedman~\cite{nF83}.
\end{paragr}

\begin{paragr}\textbf{(Strong Mixing).}
A transformation $T$ is \textbf{mixing}
when for any measurable sets $A$ and $B$,
\[
\lim_{n\to\infty}\mu(T^{n}(A) \cap B) - \mu(A)\mu(B) = 0.
\]
Given a transformation $T$, a strictly increasing sequence 
of positive integers $\{t_{n}\}$ is \textbf{mixing with respect 
to $\boldsymbol{T}$} when for any measurable sets $A$ and $B$,
\[
\lim_{n\to\infty}\mu(T^{t_{n}}(A) \cap B) - \mu(A)\mu(B) = 0.
\]
Mixing and ergodicity on sequences are related by the Blum-Hanson 
theorem \cite{BH60}.
\begin{thm}[Blum-Hanson]
Let $T$ be an ergodic transformation.  Then $T$ is a mixing 
transformation if and only if every strictly increasing sequence of 
integers is ergodic with respect to $T$.
\end{thm}
\end{paragr}

\begin{paragr}\textbf{Weak Mixing.}
A transformation $T$ is \textbf{weak mixing} when there exists a 
mixing sequence with respect to $T$.  The following conditions are equivalent to 
weak mixing: i) there exists a density one mixing sequence with 
respect to $T$; 
ii) for any sequence 
of measurable sets $\{A_{n}\}$ and any measurable set $B$, 
\[
\lim_{n\to\infty}\frac{1}{n}\sum_{i=0}^{n-1}\big{|}\mu(T^{i}(A_{n}) \cap B) - 
\mu(A_{n})\mu(B)\big{|} = 0;
\]
iii) the 
transformation $T\times T$ is ergodic; and iv) \textbf{double ergodicity}, 
for any measurable sets $A$ and $B$ of positive measure there exists a positive integer 
$n$ such that $\mu(T^{n}(A) \cap A)\mu(T^{n}(A) \cap B)>0$.  
The first three were known to von Neumann and Kakutani, the last was 
first proved in \cite{hF81} and later in \cite{BFMS01} where it was shown 
for staircase transformations.
Given an ergodic transformation $T$, a sequence of integers 
$\{a_{n}\}$ is \textbf{weak mixing with respect to 
$\boldsymbol{T}$} when for any sequence of measurable sets 
$\{A_{n}\}$ and any measurable set $B$,
\[
\lim_{n\to\infty}\frac{1}{n}\sum_{i=0}^{n-1}\big{|}\mu(T^{a_{i}}(A_{n}) 
\cap B) - \mu(A_{n})\mu(B)\big{|} = 0.
\]
A consequence of our generalization in section~\ref{S:dynavg}
(Proposition~\ref{P:wmimperg}) is that any weak mixing sequence (with respect 
to $T$) must be ergodic (with respect to $T$).
\end{paragr}

\begin{paragr}\textbf{Notation.}
To show that expressions 
converge to the same limit, we use the $\pm$ notation in the 
following manner.  The reader may verify that the $\pm$ ``quantity'' 
obeys the ``rules'' of arithmetic.
\begin{notation}
For any quantities $\boldsymbol{A}$, $\boldsymbol{B}$ and 
$\boldsymbol{C}$, the notation
\[
\boldsymbol{A} = \boldsymbol{B} \pm 
\big{(}\boldsymbol{C}\big{)}\quad\quad\text{shall mean}\quad\quad
\big{|}\boldsymbol{A} - \boldsymbol{B}\big{|} \leq \big{|}\boldsymbol{C}\big{|}.
\]
\end{notation}
\end{paragr}

\section{Dynamical Sequences}

\begin{paragr}\textbf{Basic Notions.}
We introduce the notion of dynamical sequences of 
integers which arises when considering the spacer 
levels added to each column in the construction of rank one 
transformations.  We use dynamical sequences to generalize the Ces\`{a}ro 
averaging over a static sequence in the definitions of the previous section to 
a moving Ces\`{a}ro average.
\begin{definition}\label{D:dynseq}
Given a nondecreasing sequence of positive integers
$\{r_{n}\}$, a \textbf{\emph{dynamical sequence}} of integers
$\{s_{n,i}\}_{\{r_{n}\}}$ is an indexed collection of 
integers $s_{n,i}$ for $n \geq 0$ and $0 \leq i < r_{n}$.
\end{definition}
\begin{definition}
Let $\{s_{n,i}\}_{\{r_{n}\}}$ be a dynamical sequence of integers.  
The \textbf{\emph{sequence of averages}} for 
$\{s_{n,i}\}_{\{r_{n}\}}$, denoted $\{\bar{s}_{n}\}$, and 
the \textbf{\emph{sequence of ranges}} for $\{s_{n,i}\}_{\{r_{n}\}}$, denoted 
$\{\widetilde{s}_{n}\}$, are given by 
\[
\bar{s}_{n} = 
\Big{\lfloor}\frac{1}{r_{n}}\sum_{i=0}^{r_{n}-1}s_{n,i}\Big{\rfloor}
\quad\text{and}\quad
\widetilde{s}_{n} = \max_{0\leq 
i<r_{n}}s_{n,i} - \min_{0\leq i<r_{n}}s_{n,i}
\]
for all positive integers $n$.
The \textbf{\emph{representative dynamical sequence}} of 
$\{s_{n,i}\}_{\{r_{n}\}}$, denoted $\{\widehat{s}_{n,i}\}_{\{r_{n}\}}$ is 
given by $\widehat{s}_{n,i} = s_{n,i} - \bar{s}_{n}$ for all $n$ and $i$.
\end{definition}
Note that the sequence of ranges for the representative dynamical 
sequence of a given dynamical sequence is equal to that of the 
dynamical sequence.
\begin{definition}
Given a dynamical sequence of integers $\{s_{n,i}\}_{\{r_{n}\}}$, the \textbf{\emph{
family of partial sum dynamical sequences}}, denoted 
$\{s_{n,i}^{(k)}\}_{\{r_{n}^{(k)}\}}$, is given by $r_{n}^{(k)} = r_{n} - k$ and 
$s_{n,i}^{(k)}=\sum_{z=0}^{k-1}s_{n,i+z}$ for all integers $n$, $i$ and $k$ 
such that $0 \leq i < i + k < r_{n}$.
\end{definition}
\begin{definition}
A dynamical sequence of integers 
$\{s_{n,i}\}_{\{r_{n}\}}$ is \textbf{\emph{pathological}} when the 
sequence $\{r_{n}\}$ has a finite limit point 
($\liminf_{n\to\infty}r_{n} < \infty$).
\end{definition}
We assume for the remainder that all dynamical sequences are not 
pathological (see Proposition~\ref{P:rnbounded}).
\end{paragr}

\begin{paragr}\textbf{Monotonic Dynamical Sequences.}
We generalize the property of monotonicity used in the Blum-Hanson theorem 
to dynamical sequences yielding three distinct, related concepts.  
(The symbol $\#$ denotes cardinality.)
\begin{definition}
A dynamical sequence of integers 
$\{s_{n,i}\}_{\{r_{n}\}}$ is \textbf{\emph{strictly increasing}} when for all 
$n$ and all $0\leq i<r_{n}-1$, the term $s_{n,i}< s_{n,i+1}$.  The 
dynamical sequence is \textbf{\emph{nondecreasing}} when for all $n$ 
and all $0\leq i < r_{n}-1$, the term $s_{n,i} \leq s_{n,i+1}$.
\end{definition}
\begin{definition}
A dynamical sequence of integers $\{s_{n,i}\}_{\{r_{n}\}}$ is 
\textbf{\emph{square monotone}} when for any fixed positive integer 
$M$,
\[
\lim_{n\to\infty}\frac{1}{r_{n}^{2}}\#\{(i,j)\in\mathbb{Z}_{r_{n}}\times\mathbb{Z}_{r_{n}} : \big{|} s_{n,i} - 
s_{n,j}\big{|} < M\} = 0.
\]
\end{definition}
\begin{definition}
A dynamical sequence of integers $\{s_{n,i}\}_{\{r_{n}\}}$ is 
\textbf{\emph{weak monotone}} when for any fixed positive integer 
$M$, 
\[
\lim_{n\to\infty}\frac{1}{r_{n}}\#\{i\in\mathbb{Z}_{r_{n}} :
\big{|} s_{n,i} \big{|} < M\} = 0.
\]
\end{definition}
\begin{proposition}\label{P:incr}
If a given dynamical sequence is strictly increasing, then it is 
square monotone; if it is square monotone, then it is weak monotone.
\end{proposition}
\begin{proof}
Given a strictly increasing dynamical 
sequence of integers $\{s_{n,i}\}_{\{r_{n}\}}$,
for any $n$ and any $0\leq j\leq 
i < r_{n}$, we have that
$s_{n,i} - s_{n,j} > s_{n,i-1} - s_{n,j} \geq 1 + s_{n,i-1} - 
s_{n,j} > \ldots \geq (i - j) + s_{n,i-(i-j)} - s_{n,j} = i - j$;
thus, $\{s_{n,i}\}_{\{r_{n}\}}$ is square monotone.  That a 
square monotone dynamical sequence is weak monotone is left to the reader.
\end{proof}
\end{paragr}

\begin{paragr}\textbf{Dynamical Subsequences and Multiplicity.}
Let $\{s_{n,i}\}_{\{r_{n}\}}$ be a dynamical sequence of integers.  
The \textbf{multiplicity functions} of 
$\{s_{n,i}\}_{\{r_{n}\}}$, denoted $\{R_{n}\}$, are given 
by $R_{n}(\ell) = \#\{i\in\mathbb{Z}_{r_{n}}:s_{n,i} = \ell\}$ for 
each integer $\ell$.
If a dynamical sequence of integers $\{a_{n,i}\}_{\{q_{n}\}}$ with 
multiplicity functions $\{R_{n}^{\prime}\}$ has the 
property that for any integers $n$ and $\ell$, $R_{n}^{\prime}(\ell) \leq 
R_{n}(\ell)$, then $\{a_{n,i}\}_{\{q_{n}\}}$ is a 
\textbf{dynamical subsequence} of $\{s_{n,i}\}_{\{r_{n}\}}$.  
The multiplicity functions for the partial sum dynamical 
sequences are handled similarly and denoted by $\{R_{n}^{(k)}\}$ 
for each positive integer $k$.
\end{paragr}

\begin{paragr}\textbf{Density of Dynamical Sequences.}
For the following definitions replace $\lim$ by $\limsup$ or 
$\liminf$ for notions of upper and lower density in the case when the 
limit does not exist.
\begin{definition}
A dynamical sequence $\{s_{n,i}\}_{\{r_{n}\}}$ has 
\textbf{\emph{density}} given by $D = 
\lim_{n\to\infty}\frac{r_{n}}{\widetilde{s}_{n}}$
where $\{\widetilde{s}_{n}\}$ is the sequence of ranges for 
$\{s_{n,i}\}_{\{r_{n}\}}$ given by $\widetilde{s}_{n} = \max_{i}s_{n,i} - 
\min_{i}s_{n,i}$.
\end{definition}
In particular, we refer to dynamical sequences being of positive (lower) 
density and of finite (upper) density.
\begin{definition}
A dynamical sequence of integers $\{s_{n,i}\}_{\{r_{n}\}}$ has 
\textbf{\emph{density in $\boldsymbol{\mathbb{Z}}$}} given by
$D_{\mathbb{Z}} = 
\lim_{n\to\infty}\frac{1}{\widetilde{s}_{n}}\#\{\ell\in\mathbb{Z}:s_{n,i} = 
\ell \quad \text{for some $i\in\mathbb{Z}_{r_{n}}$}\}$.
\end{definition}
Note that $D_{\mathbb{Z}} \leq D$ for all any dynamical sequence and that
when the $s_{n,i}$ take on distinct values the density of 
the sequence equals the density in $\mathbb{Z}$.
\begin{definition}
A dynamical subsequence $\{a_{n,i}\}_{\{q_{n}\}}$ of a dynamical 
sequence $\{s_{n,i}\}_{\{r_{n}\}}$ has 
\textbf{\emph{density in $\{s_{n,i}\}_{\{r_{n}\}}$}} 
given by $D_{s} = 
\lim_{n\to\infty}\frac{q_{n}}{r_{n}}$.
\end{definition}
\end{paragr}

\section{Dynamical Ces\`{a}ro Averaging}\label{S:dynavg}

\begin{paragr}\textbf{Dynamical Sequence Ergodicity.}
The application of dynamical sequences to expressions like those found 
in the ergodic theorems and the definition of ergodic sequence is 
crucial to our main result.  We explore notions of ergodicity and weak 
mixing for dynamical sequences and present a generalized version of the 
Blum-Hanson theorem and a density result for weak mixing 
transformations.
\begin{definition}
Given an ergodic transformation $T$, 
a dynamical sequence of nonnegative integers 
$\{s_{n,i}\}_{\{r_{n}\}}$ is 
\textbf{\emph{ergodic with respect to $\boldsymbol{T}$}} when for any 
measurable set $B$,
\[
\lim_{n\to\infty} \int\big{|}\frac{1}{r_{n}}\sum_{i=0}^{r_{n}-1}\chi_{B}\circ 
T^{-s_{n,i}} - \mu(B)\big{|}^{2} d\mu = 0.
\]
We refer to the above sequence of functions as \textbf{\emph{dynamical Ces\`{a}ro 
averages}} for the $L^{2}$ function $\chi_{B}$ over the dynamical sequence $\{s_{n,i}\}_{\{r_{n}\}}$.
\end{definition}
\begin{definition}
Given an ergodic transformation $T$, a dynamical sequence of integers
$\{s_{n,i}\}_{\{r_{n}\}}$
is \textbf{\emph{uniformly ergodic with respect to $\boldsymbol{T}$}} when 
for any sequence of positive integers 
$\{k_{n}\}$ such that $\limsup_{n\to\infty}\frac{k_{n}}{r_{n}} < 1$ 
and any measurable set $B$,
\[
\lim_{n\to\infty} \int 
\big{|}\frac{1}{r_{n}^{(k_{n})}}\sum_{i=0}^{r_{n}^{(k_{n})}-1}\chi_{B}\circ 
T^{-s_{n,i}^{(k_{n})}} - \mu(B)\big{|}^{2} d\mu = 0.
\]
Equivalently, the above limit converges to zero uniformly over all 
positive integers $k$ such that $\frac{k}{r_{n}}$ is bounded below 
one---given any $\epsilon >0$ there exists an integer $N$ such that 
for all integers $n \geq N$ and all integers $k$ such that 
$\limsup_{n\to\infty}\frac{k}{r_{n}}<1$, the integral above is less 
than $\epsilon$.
\end{definition}
Note that uniform ergodicity implies ergodicity for a given dynamical sequence
with respect to a given transformation.
\begin{definition}
Given an ergodic transformation $T$, a 
dynamical sequence of integers
$\{s_{n,i}\}_{\{r_{n}\}}$
is \textbf{\emph{weak ergodic with respect to $\boldsymbol{T}$}} when
for any measurable sets $A$ and $B$,
\[
\lim_{n\to\infty} \frac{1}{r_{n}}\sum_{i=0}^{r_{n}-1}
\mu(T^{s_{n,i}}(A) \cap B) - \mu(A)\mu(B) = 0.
\]
\end{definition}
Note that a dynamical sequence that is ergodic with respect to a given 
transformation must be weak ergodic with respect to the transformation.
\end{paragr}

\begin{paragr}\textbf{Weak Mixing Dynamical Sequences.}
Similar to the notions of ergodicity and uniform ergodicity on 
dynamical sequences, we 
define weak mixing and uniform weak mixing with respect to a given 
transformation.
\begin{definition}
Given an ergodic transformation $T$, a 
dynamical sequence of integers
$\{s_{n,i}\}_{\{r_{n}\}}$
is \textbf{\emph{weak mixing with respect to $\boldsymbol{T}$}} when
for any sequence of measurable sets $\{A_{n}\}$ and any measurable 
set $B$,
\[
\lim_{n\to\infty}\frac{1}{r_{n}}\sum_{i=0}^{r_{n}-1}
\big{|}\mu(T^{s_{n,i}}(A_{n}) \cap B) - \mu(A_{n})\mu(B)\big{|} = 0.
\]
\end{definition}
\begin{definition}
Given an ergodic transformation $T$, a dynamical sequence of integers
$\{s_{n,i}\}_{\{r_{n}\}}$
is \textbf{\emph{uniformly weak mixing with respect to $\boldsymbol{T}$}} when 
for any sequence of measurable sets $\{A_{n}\}$, and measurable set $B$ and any sequence of positive integers 
$\{k_{n}\}$ such that $\limsup_{n\to\infty}\frac{k_{n}}{r_{n}} < 1$,
\[
\lim_{n\to\infty}\frac{1}{r_{n}^{(k_{n})}}\sum_{i=0}^{r_{n}^{(k_{n})}-1}
\big{|}\mu(T^{s_{n,i}^{(k_{n})}}(A_{n}) \cap B) - \mu(A_{n})\mu(B)\big{|} = 0.
\]
Equivalently, the limit above converges uniformly over all measurable 
sets $A$ and all positive integers $k$ such that $\frac{k}{r_{n}}$ is 
bounded below one.
\end{definition}
\begin{proposition}\label{P:wmimperg}
Let $\{s_{n,i}\}_{\{r_{n}\}}$ be a dynamical sequence that is weak mixing with 
respect to an ergodic transformation $T$.  Then 
$\{s_{n,i}\}_{\{r_{n}\}}$ is ergodic with respect to $T$.
\end{proposition}
\begin{proof}
Let $\{s_{n,i}\}_{\{r_{n}\}}$ and $T$ be as above.  Fix a measurable 
set $B$.  For any positive integer $n$, let $A_{n}^{(+)} = 
\{x : \frac{1}{r_{n}}\sum_{i=0}^{r_{n}-1}\chi_{B}\circ 
T^{-s_{n,i}}(x) - \mu(B) \geq 0\}$ and $A_{n}^{(-)}$ similarly.  Then,
{\allowdisplaybreaks
\begin{align*}
\int\big{|}\frac{1}{r_{n}}&\sum_{i=0}^{r_{n}-1}\chi_{B}\circ 
T^{-s_{n,i}} - \mu(B)\big{|} d\mu \\
&= \big{|}\int_{A_{n}^{(+)}}\frac{1}{r_{n}}\sum_{i=0}^{r_{n}-1}\chi_{B}\circ 
T^{-s_{n,i}} - \mu(B) d\mu\big{|} \\ &\quad\quad\quad + \big{|}\int_{A_{n}^{(-)}}\frac{1}{r_{n}}\sum_{i=0}^{r_{n}-1}\chi_{B}\circ 
T^{-s_{n,i}} - \mu(B) d\mu\big{|} \\
&\leq 
\frac{1}{r_{n}}\sum_{i=0}^{r_{n}-1}\big{|}\mu(T^{s_{n,i}}(A_{n}^{(+)}) \cap 
B)\big{|} + \frac{1}{r_{n}}\sum_{i=0}^{r_{n}-1}\big{|}\mu(T^{s_{n,i}}(A_{n}^{(-)}) \cap 
B)\big{|}
\end{align*}}
which approaches zero as $n\to\infty$ since $\{s_{n,i}\}_{\{r_{n}\}}$ 
is weak mixing with respect to $T$.
\end{proof}
\end{paragr}

\begin{paragr}\textbf{Weak Mixing and Density on Dynamical Sequences.}
Furstenberg's results relating sequences of positive density and weak 
mixing \cite{hF81} are generalized to dynamical sequences as follows.
\begin{theorem}\label{T:weakmixdens}
Let $\{a_{n,i}\}_{\{q_{n}\}}$ be a dynamical subsequence of a dynamical sequence 
of integers $\{s_{n,i}\}_{\{r_{n}\}}$ that has positive (lower) density 
in $\{s_{n,i}\}_{\{r_{n}\}}$.  If $\{s_{n,i}\}_{\{r_{n}\}}$ is weak 
mixing with respect to an ergodic transformation $T$, then 
$\{a_{n,i}\}_{\{q_{n}\}}$ is weak mixing with respect to $T$.
\end{theorem}
\begin{proof}
Let $\{a_{n,i}\}_{\{q_{n}}$, $\{s_{n,i}\}_{\{r_{n}\}}$ and $T$ be as 
above and let $D = \liminf_{n\to\infty}\frac{q_{n}}{r_{n}}$ be the 
density.  Note that $\limsup_{n\to\infty}\frac{r_{n}}{q_{n}} = 
\frac{1}{D}$.  Let $\{A_{n}\}$ be any sequence of measurable sets and 
$B$ any measurable set.  Then
\begin{align*}
\frac{1}{q_{n}}\sum_{i=0}^{q_{n}-1}\big{|}&\mu(T^{a_{n,i}}(A_{n}) \cap 
B) - \mu(A_{n})\mu(B)\big{|} \\ &\leq 
\frac{1}{q_{n}}\sum_{i=0}^{r_{n}-1}\big{|}
\mu(T^{s_{n,i}}(A_{n}) \cap B) - \mu(A_{n})\mu(B)\big{|} \\ &< 
\frac{1}{D}\frac{1}{r_{n}}\sum_{i=0}^{r_{n}-1}\big{|}\mu(T^{s_{n,i}}(A_{n}) \cap B) - \mu(A_{n})\mu(B)\big{|}
\end{align*}
approaches zero as $n\to\infty$ if $\{s_{n,i}\}_{\{r_{n}\}}$ is weak 
mixing with respect to $T$.
\end{proof}
\begin{corollary}\label{C:pdimpwm}
Let $T$ be a weak mixing transformation and let 
$\{s_{n,i}\}_{\{r_{n}\}}$ be a dynamical sequence of integers.  If 
$\{s_{n,i}\}_{\{r_{n}\}}$ has positive (lower) density and takes on 
each value $\ell$ no more than once for each $n$, then 
$\{s_{n,i}\}_{\{r_{n}\}}$ is weak mixing with respect to $T$.
\end{corollary}

\paragraph{Generalized Blum-Hanson Theorem}
The following pair of theorems, generalizations of the Blum-Hanson 
theorem and its corresponding weak version,  
characterize the ergodicity of dynamical sequences with respect to mixing 
transformations.
\begin{theorem}\label{T:genBH}
Let $T$ be an ergodic transformation.  Then $T$ is mixing if and only 
if every square monotone dynamical sequence of integers
is ergodic with respect to $T$.
\end{theorem}
\begin{proof}
Let $T$ be a mixing transformation and $\{s_{n,i}\}_{\{r_{n}\}}$ a 
square monotone dynamical sequence of integers.  For any measurable set $B$ and 
any $\epsilon > 0$ there exists $M > 0$ such that for all $m\geq M$ 
or $m \leq -M$, 
$\big{|}\mu(T^{m}(B)\cap B)-\mu(B)\mu(B)\big{|} < \epsilon$.  Since 
$\{s_{n,i}\}_{\{r_{n}\}}$ is square monotone, there exists $N > 
0$ such that for all $n\geq N$, 
$\#\{(i,j)\in\mathbb{Z}_{r_{n}}\times\mathbb{Z}_{r_{n}}:\big{|} s_{n,i}-s_{n,j}\big{|}<M\} < 
\epsilon r_{n}^{2}$.  Then,
{\allowdisplaybreaks
\begin{align*}
\int&\big{|}\frac{1}{r_{n}}\sum_{i=0}^{r_{n}-1}\chi_{B}\circ T^{-s_{n,i}} - 
\mu(B)\big{|}^{2} d\mu \\ &= 
\frac{1}{r_{n}^{2}}\sum_{i,j=0}^{r_{n}-1}\mu(T^{s_{n,i}}(B)\cap 
T^{s_{n,j}}(B)) - \mu(B)\mu(B) \\
&\leq 
\frac{1}{r_{n}^{2}}\sum_{i,j=0}^{r_{n}-1}\mu(T^{s_{n,i}-s_{n,j}}(B)\cap 
B) - \mu(B)\mu(B) \\
&\leq \frac{1}{r_{n}^{2}}\sum_{|s_{n,i}-s_{n,j}|\geq 
M}\big{|}\mu(T^{s_{n,i}-s_{n,j}}(B)\cap B) - \mu(B)\mu(B)\big{|} + 
\frac{1}{r_{n}^{2}}\epsilon r_{n}^{2}\mu(B) \\
&< \frac{1}{r_{n}^{2}}\sum_{|s_{n,i}-s_{n,j}|\geq M}\epsilon + 
\epsilon\mu(B) \leq \epsilon (1+\mu(B))
\end{align*}}
which approaches zero by letting $\epsilon\to 0$.  Hence, 
$\{s_{n,i}\}_{\{r_{n}\}}$ is ergodic with respect to $T$.

Since mixing is equivalent to R\'{e}nyi mixing: for all measurable sets $B$,
$\mu(T^{m}(B)\cap B)\to\mu(B)\mu(B)$, if $T$ is not mixing then 
there exists a measurable set $B$, a $\delta > 0$ and a strictly 
increasing sequence of positive integers $\{t_{m}\}$ such that for 
all $m$, $\mu(T^{t_{m}}(B)\cap B) - \mu(B)\mu(B) \geq \delta$.  
Define the dynamical sequence of integers $\{s_{n,i}\}_{\{r_{n}\}}$ by $r_{n} = n$ 
and $s_{n,i} = t_{i}$.  Then, first using the H\"{o}lder Inequality 
and then the triangle inequality twice,
{\allowdisplaybreaks
\begin{align*}
\int\big{|}\frac{1}{r_{n}}\sum_{i=0}^{r_{n}-1}\chi_{B}\circ 
T^{-s_{n,i}} - \mu(B)\big{|}^{2} d\mu &\geq \int\big{|}\frac{1}{r_{n}}\sum_{i=0}^{r_{n}-1}\chi_{B}\circ 
T^{-s_{n,i}} - \mu(B)\big{|}d\mu \\
&\geq \int_{B}\big{|}\frac{1}{n}\sum_{i=0}^{n-1}\chi_{B}\circ 
T^{-t_{i}} - \mu(B)\big{|}d\mu \\
&\geq \frac{1}{n}\sum_{i=0}^{n-1}\mu(T^{t_{i}}(B)\cap B) - 
\mu(B)\mu(B) \geq \delta.
\end{align*}}
Thus, $\{s_{n,i}\}_{\{r_{n}\}}$ is not ergodic with respect to $T$.  
It remains only to show that $\{s_{n,i}\}_{\{r_{n}\}}$ is weakly 
increasing.  But, since $\{t_{m}\}$ is strictly increasing, 
$\{s_{n,i}\}_{\{r_{n}\}}$ is strictly increasing and thus square 
monotone by Proposition~\ref{P:incr}.
\end{proof}
\begin{theorem}\label{T:genweakBH}
Let $T$ be an ergodic transformation.  Then $T$ is mixing if and only 
if every weak monotone dynamical sequence of integers
is weak ergodic with respect to $T$.
\end{theorem}
\begin{proof}
Let $T$ be a mixing transformation and $\{s_{n,i}\}_{\{r_{n}\}}$ a 
weak monotone dynamical sequence.  For any measurable sets $A$ and $B$ 
and any $\epsilon > 0$ there exists a positive integer $M$ such that 
for all integers $m \geq M$ or $m \leq -M$, $\big{|}\mu(T^{m}(A) \cap B) - 
\mu(A)\mu(B)\big{|} < \epsilon$.  There also exists a positive integer $N$ 
such that for all integers $n \geq N$, 
$\#\{i\in\mathbb{Z}_{r_{n}}:\big{|}s_{n,i}\big{|}<M\} < \epsilon r_{n}$.  Hence,

{\allowdisplaybreaks
\begin{align*}
\big{|}\frac{1}{r_{n}}&\sum_{i=0}^{r_{n}-1}\mu(T^{s_{n,i}}(A) \cap B) - 
\mu(A)\mu(B)\big{|} \\
&\leq \big{|}\frac{1}{r_{n}}\sum_{\substack{i=0 \\ 
|s_{n,i}|<M}}^{r_{n}-1}\mu(T^{s_{n,i}}(A) \cap B) - \mu(A)\mu(B)\big{|} 
\\ &\quad\quad\quad + 
\big{|}\frac{1}{r_{n}}\sum_{\substack{i=0 \\ |s_{n,i}|\geq 
M}}^{r_{n}-1}\mu(T^{s_{n,i}}(A) \cap B) - \mu(A)\mu(B)\big{|} \\
&< \frac{1}{r_{n}}\epsilon r_{n}\mu(A) + 
\frac{1}{r_{n}}r_{n}\epsilon = \epsilon (1+\mu(A))
\end{align*}}
which approaches zero by letting $\epsilon \to 0$.  Hence, 
$\{s_{n,i}\}_{\{r_{n}\}}$ is weak ergodic with respect to $T$.
Conversely, if $T$ is not mixing then there exists a strictly 
increasing sequence $\{t_{n}\}$, $\delta > 0$ and measurable sets $A$ 
and $B$ such that $\mu(T^{t_{n}}(A) \cap B) - \mu(A)\mu(B) \geq 
\delta$ for all $n$.  For all $n$ and $i$, set $r_{n} = n$ and 
$s_{n,i} = t_{i}$.  Then $\{s_{n,i}\}_{\{r_{n}\}}$ is strictly 
increasing and so is weak monotone by Proposition~\ref{P:incr}.  Clearly, 
$\{s_{n,i}\}_{\{r_{n}\}}$ is not weak ergodic with respect to $T$.  
\end{proof}
\end{paragr}

\section{Rank One Transformations}\label{S:rankone}

\begin{paragr}\textbf{Construction of Rank One Transformations.}
Rank one transformations are a class of ergodic transformations on 
the unit interval in $\mathbb{R}$ under standard Lebesgue measure constructed as follows.
An ordered collection of intervals all the same length is termed a 
\textbf{column} and each interval a \textbf{level} 
where the \textbf{height} of the column is the number of 
levels in the column.  The associated \textbf{column map} is 
defined by mapping each interval to the interval 
\textbf{above} (next in the order on the collection of 
intervals) it, hence the column map is defined from all but the 
\textbf{top} (last in the order) level onto all but the 
\textbf{bottom} (first in the order) level.

We describe the procedure for \textbf{cutting and stacking} a 
column, $C_{n}$, to obtain a new column, $C_{n+1}$.
Fix a column $C_{n}$ with height $h_{n}$, levels $I_{n,j}$ and column 
map $T_{n}$:
$C_{n} = \{I_{n,j}\}_{j=0}^{h_{n}-1}$ where $T_{n}(I_{n,j}) = I_{n,j+1}$
for $j\ne h_{n}-1$.
For some given $r_{n} > 0$, \textbf{cut} $C_{n}$ into $r_{n}$ subcolumns by 
cutting each level $I_{n,j}$ into $r_{n}$ \textbf{sublevels}, 
$\{I_{n,j}^{[i]}\}_{i=0}^{r_{n}-1}$, of equal length, 
$\frac{1}{r_{n}}\mu(I_{n,j})$, where 
$I_{n,j}^{[0]}$ is the leftmost sublevel and $I_{n,j}^{[r_{n}-1]}$ is 
the rightmost.  Then, the \textbf{subcolumns} of $C_{n}$ are $C_{n}^{[i]} = 
\{I_{n,j}^{[i]}\}_{j=0}^{h_{n}-1}$.  By 
preserving the order on the levels, each subcolumn, $C_{n}^{[i]}$, is 
a column in its own right with the associated map $T_{n}^{[i]}$ which is the 
restriction of $T_{n}$ to $C_{n}^{[i]}$.

Given an indexed collection of nonnegative integers $\{s_{n,i}\}_{i=0}^{r_{n}-1}$, 
the spacer values for $C_{n}$, place \textbf{spacer levels} (``new'' 
intervals the size of each sublevel) above each subcolumn 
by adding $s_{n,i}$ levels above $C_{n}^{[i]}$ and 
\textbf{stack} the resulting subcolumns with spacers right on top of left 
yielding a new column $C_{n+1}$ with height $h_{n+1} = r_{n}h_{n} + 
\sum_{i=0}^{r_{n}-1}s_{n,i}$.  Denote the \textbf{union of spacer 
levels} added to $C_{n}$ by $S_{n}$, the collection of levels in 
$C_{n+1}$ that are not sublevels of levels in $C_{n}$. The associated column map $T_{n+1}$ 
restricts to $T_{n}$ on the levels in $C_{n}$, as above, and extends 
it to the spacer levels added as well as all but the leftmost sublevel 
of the bottom level and the topmost spacer level over the rightmost 
subcolumn of $C_{n}$.  We will use implicitly the following facts in the sequel.
\begin{lemma}
For any sublevel $I_{n,j}^{[i]}$ in $C_{n}$, $I_{n,j}^{[i]} = 
I_{n+1,j+ih_{n}+\sum_{z=0}^{i-1}s_{n,z}}$ is a level
in $C_{n+1}$.
\end{lemma}
\begin{lemma}
For any sublevel $I_{n,j}^{[i]}$ in $C_{n}$ where $i \ne r_{n}-1$,
$T_{n+1}^{h_{n}+s_{n,i}}(I_{n,j}^{[i]}) = I_{n,j}^{[i+1]}$.
\end{lemma}
Thus, given a nondecreasing sequence of positive integers $\{r_{n}\}$, the \textbf{sequence 
of cuts}, and a dynamical sequence of nonnegative 
integers $\{s_{n,i}\}_{\{r_{n}\}}$, the \textbf{sequence of 
spacers}, we may 
construct an infinite sequence of columns $\{C_{n}\}$ where $C_{0} = 
\{I_{0,0}\}$ is defined to be a single level (whose length will be 
chosen later to normalize the transformation) and each $C_{n+1}$ is 
constructed by cutting and stacking each $C_{n}$ as described above 
using $r_{n}$ cuts and $\{s_{n,i}\}_{r_{n}}$ spacers.  The 
\textbf{sequence of heights} is then defined recursively by $h_{0} = 
1$ and $h_{n+1} = r_{n}h_{n} + \sum_{i=0}^{r_{n}-1}s_{n,i}$.

Since the initial level $I_{0,0}$ has some finite length, the length 
of each level in $C_{n}$ approaches zero as $n$ becomes large.  
Hence, the associated sequence of column maps, $\{T_{n}\}$, approaches a map $T$ defined on all but a 
measure zero subset of the union of the initial levels and all the 
spacer levels added at each column.  Note that $T_{n}$ is the 
restriction of $T$ to the levels in $C_{n}$ (except the top) for each 
$n\geq 0$.  When $T$ is defined on a finite 
measure space, choose $I_{0,0}$ to have length such that $T$ is defined 
on the unit interval.

A transformation $T$ is formally a \textbf{rank one transformation} 
when $T$ can be realized as the limit of a sequence of maps defined 
by the cut and stack construction applied repeatedly to a single column as described above.
The reader may verify that rank one transformations are invertible, 
measurable, measure-preserving and ergodic since any measurable set is contained in a 
union of levels.  
Given a rank one transformation $T$ with cut sequence $\{r_{n}\}$ and 
spacer sequence $\{s_{n,i}\}_{\{r_{n}\}}$, we define 
the \textbf{sequence of spacer averages} for $T$ to be the sequence of 
averages for $\{s_{n,i}\}_{\{r_{n}\}}$, denoted $\{\bar{s}_{n}\}$, given by $\bar{s}_{n} = 
\big{\lfloor}\frac{1}{r_{n}}\sum_{i=0}^{r_{n}-1}s_{n,i}\big{\rfloor}$ 
and the 
\textbf{sequence of window heights} for $T$ to be the sequence of 
positive integers $\{w_{n}\}$ given by $w_{n} = h_{n} + \bar{s}_{n}$.  Note that 
$\lim_{n\to\infty}\frac{w_{n}}{h_{n}} = 1$ when $T$ is finite 
measure-preserving.
We define the \textbf{representative spacer sequence} for $T$ to be the 
representative dynamical sequence 
of $\{s_{n,i}\}_{\{r_{n}\}}$, denoted $\{\widehat{s}_{n,i}\}_{\{r_{n}\}}$, given by $s_{n,i} = 
\widehat{s}_{n,i} - \bar{s}_{n}$.  
Then the average value of the representative spacer sequence for each $n$ is 
$\frac{1}{r_{n}}\sum_{i=0}^{r_{n}-1}\widehat{s}_{n,i} = 
\frac{1}{r_{n}}\sum_{i=0}^{r_{n}-1}s_{n,i} - 
\bar{s}_{n}$ is between zero and one.
\begin{lemma}
For any sublevel $I_{n,j}^{[i]}$ in $C_{n}$ where $i \ne r_{n}-1$,
$T^{w_{n}}(I_{n,j}^{[i]}) = T^{-\widehat{s}_{n,i}}(I_{n,j}^{[i+1]})$.
\end{lemma}
\begin{definition}
A rank one transformation $T$ has an \textbf{\emph{ergodic spacer 
sequence}} when the spacer sequence for $T$ is ergodic with respect 
to $T$.  Similarly, we define \textbf{\emph{weakly ergodic spacer 
sequence}}, \textbf{\emph{uniformly ergodic spacer 
sequence}}, \textbf{\emph{weak mixing spacer sequence}} and 
\textbf{\emph{uniformly weak mixing spacer sequence}}.
\end{definition}
Note that since $T$ is measure-preserving the spacer sequence is 
ergodic (respectively, uniformly ergodic or (uniformly) weak mixing) with 
respect to $T$ if and 
only if the representative spacer sequence is ergodic (respectively, 
uniformly ergodic or (uniformly) weak mixing) with respect to $T$.
We will use implicitly the following 
observation in the sequel; the proof is standard.
\begin{proposition}\label{P:rnbounded}
Let $T$ be a rank one transformation with sequence of cuts 
$\{r_{n}\}$ having a finite limit point.  Then $T$ is partially rigid and therefore cannot be 
mixing.  Hence, if the spacer sequence for $T$ is a pathological 
dynamical sequence then $T$ is not mixing. 
\end{proposition}
\end{paragr}

\begin{paragr}\textbf{Restricted Growth Rank One Transformations.}
We introduce the class of restricted growth rank one transformations 
characterized by adding spacer levels whose maximum variation in height approaches 
zero relative to the height of the column being cut and stacked.  Note 
that adding spacer levels whose total height approaches zero 
relative to the height of the column \emph{resulting from the cut and 
stack procedure} is a necessary condition for the space
the transformation is defined on to be finite.
\begin{definition}
Given a rank one transformation $T$  with representative spacer sequence 
$\{\widehat{s}_{n,i}\}_{\{r_{n}\}}$ and height sequence $\{h_{n}\}$, the 
transformation $T$ has
\textbf{\emph{restricted growth}}
when for any sequences of positive integers $\{i_{n}\}$ and 
$\{k_{n}\}$ such that $0\leq i_{n} < r_{n}$ and 
$\limsup_{n\to\infty}\frac{k_{n}}{r_{n}} < 1$, we have 
$\lim_{n\to\infty}\frac{1}{h_{n}}\widehat{s}_{n,i_{n}}^{(k_{n})} = 0$ where $\{\widehat{s}_{n,i}^{(k)}\}_{\{r_{n}^{(k)}\}}$ denotes the family 
of partial sum dynamical sequences for $\{\widehat{s}_{n,i}\}_{\{r_{n}\}}$; 
equivalently, $\frac{1}{h_{n}}\widehat{s}_{n,i}^{(k)} \to 0$ as $n\to\infty$ uniformly 
over positive integers $i$ and $k$ such that $0 \leq i < r_{n} - k$ 
and $\frac{k}{r_{n}}$ is bounded below one.
\end{definition}
For completeness, we provide an example of a rank one transformation 
on a finite space that does not have restricted growth.
Construct $T$ using the sequence of cuts $\{r_{n}\}$ given by $r_{n} = 
2(2^{n}-1)$ and the spacer sequence 
$\{s_{n,i}\}_{\{r_{n}\}}$ given by $s_{n,0} = 
h_{n}$ and $s_{n,i} = 0$ for $0 < i < r_{n}$ where $\{h_{n}\}$ is the 
sequence of heights.  Letting 
$\{\widehat{s}_{n,i}\}_{\{r_{n}\}}$ denote the representative spacer sequence for $T$, observe that 
$\frac{1}{h_{n}}\widehat{s}_{n,0}^{(1)} = \frac{1}{h_{n}}(h_{n} - 
\big{\lfloor}\frac{h_{n}}{r_{n}}\big{\rfloor}) = 1 - 
\frac{1}{h_{n}}\big{\lfloor}\frac{h_{n}}{r_{n}}\big{\rfloor}$ does not approach 
zero so $T$ does not have restricted growth.  That $T$ is finite 
measure-preserving is left to the reader.
\end{paragr}

\begin{paragr}\textbf{Rank One Uniform Mixing Transformations.}
We introduce the notion of uniform mixing for rank one transformations 
by considering the sums of the mixing values over increasingly fine levels.  
Our main theorem implies that a rank one transformation is mixing if and only if the 
transformation is uniform mixing, but note that the analogous result 
does not hold for sequences.
\begin{definition}
A rank one transformation $T$ is \textbf{\emph{uniformly mixing}} when for 
any measurable set $B$,
\[
\lim_{n\to\infty}\sum_{j=0}^{h_{p}-1}\big{|}\mu(T^{n}(I_{p,j}) \cap B) - 
\mu(I_{p,j})\mu(B)\big{|} = 0
\]
where $\{h_{m}\}$ is the height sequence for $T$, $p$ is the positive integer such that $h_{p}\leq n<h_{p+1}$ 
and $\{I_{p,j}\}$ are the levels in the $p$th column for $T$.
\end{definition}
\begin{proposition}\label{P:unifmix}
Let $T$ be a rank one transformation.  If $T$ is a uniformly mixing 
transformation then $T$ is a mixing transformation.
\end{proposition}
\begin{proof}
Let $T$ be a uniformly mixing rank one transformation with cut sequence $\{r_{n}\}$, 
height sequence $\{h_{n}\}$ and levels $\{I_{n,j}\}$.  For any positive integer $n$, 
let $p(n)$ denote the unique positive integer $p$ such that $h_{p} 
\leq n < h_{p+1}$.  Since the levels generate the measurable sets, it 
suffices to show that $T$ is mixing on levels.  Let $A$ and $B$ be 
unions of levels in some column $C_{N}$ for some fixed positive 
integer $N$.  Write $A = \bigcup_{j=0}^{\beta -1}I_{N,\alpha_{j}}$.  
Then, for any integer $n > h_{N+1}$,
{\allowdisplaybreaks
\begin{align*}
\big{|}\mu(&T^{n}(A) \cap B) - \mu(A)\mu(B)\big{|} \\
&= \big{|}\sum_{j=0}^{\beta -1}\mu(T^{n}(I_{N,\alpha_{j}} \cap B) - 
\mu(I_{N,\alpha_{j}})\mu(B)\big{|} \\
&\leq \sum_{j=0}^{\beta -1}\big{|}\mu(T^{n}(I_{N,\alpha_{j}}) \cap B) - 
\mu(I_{N,\alpha_{j}})\mu(B)\big{|} \\
&\leq \sum_{j=0}^{h_{N}-1}\big{|}\mu(T^{n}(I_{N,j}) \cap B) - 
\mu(I_{N,j})\mu(B)\big{|} \\
&= \sum_{j=0}^{h_{N}-1}\big{|}\sum_{z=0}^{r_{N}\ldots r_{p(n)-1} 
-1}\mu(T^{n}(I_{p(n),j + \gamma(n,z)}) \cap B) - 
\mu(I_{N,j + \gamma(n,z)})\mu(B)\big{|} \\
&\leq \sum_{j=0}^{h_{p(n)}-1}\big{|}\mu(T^{n}(I_{p(n),j}) \cap B) - 
\mu(I_{p(n),j})\mu(B)\big{|}.
\end{align*}}
\end{proof}
\begin{definition}
Given an ergodic transformation $T$, a sequence of positive 
integers $\{a_{n}\}$ is \textbf{\emph{uniformly mixing with respect to 
$\boldsymbol{T}$}} when
\[
\lim_{n\to\infty}\sum_{j=0}^{h_{p}-1}\big{|}\mu(T^{a_{n}}(I_{p,j}) \cap B) - 
\mu(I_{p,j})\mu(B)\big{|} = 0
\]
where $\{h_{m}\}$ is the height sequence for $T$, $p$ is the positive integer 
such that $h_{p}\leq a_{n}<h_{p+1}$ 
and $\{I_{p,j}\}$ are the levels in the $p$th column for $T$.
\end{definition}
Note that a uniformly mixing sequence is necessarily a mixing 
sequence as above, but that the converse statement does not hold.  The pair of 
theorems in the following section relating mixing and uniform mixing 
on height sequences to weakly ergodic and ergodic spacer sequences make 
this clear.
The concept of summing the mixing values over the levels can be 
applied to the ergodic averages as well.  The sums over the levels of 
the ergodic averages may be regarded as a specific Riemann sum for 
the ergodic integral; the following proposition makes this explicit.
\begin{proposition}\label{P:unifmixerg}
Let $T$ be a rank one transformation with height sequence $\{h_{n}\}$ 
and levels $\{I_{n,j}\}$ and let $\{s_{n,i}\}_{\{r_{n}\}}$ be a dynamical sequence of 
nonnegative integers.  Then $\{s_{n,i}\}_{\{r_{n}\}}$ is ergodic with respect 
to $T$ if 
and only if for any measurable set $B$ and any unbounded 
nondecreasing sequence of positive integers $\{p_{n}\}$ such that 
$\frac{1}{h_{p_{n}}r_{n}}\sum_{i=0}^{r_{n}-1}s_{n,i} \to 0$ as 
$n\to\infty$,
\[
\lim_{n\to\infty}\sum_{j=0}^{h_{p_{n}}-1}\big{|}\frac{1}{r_{n}}
\sum_{i=0}^{r_{n}-1}\mu(T^{-s_{n,i}}(I_{p_{n},j}) \cap B) - 
\mu(I_{p_{n},j})\mu(B)\big{|} = 0.
\]
\end{proposition}
\begin{proof}
Let $T$, $\{h_{n}\}$, $\{I_{n,j}\}$, $\{s_{n,i}\}_{\{r_{n}\}}$, $B$ 
and $\{p_{n}\}$ be 
as above.
Clearly, for each positive integer $n$, 
by the triangle inequality,
\begin{align*}
\sum_{j=0}^{h_{p_{n}}-1}\big{|}\frac{1}{r_{n}}
\sum_{i=0}^{r_{n}-1}&\mu(T^{-s_{n,i}}(I_{p_{n},j}) \cap B) - 
\mu(I_{p_{n},j})\mu(B)\big{|} \\ &\leq 
\int\big{|}\frac{1}{r_{n}}\sum_{i=0}^{r_{n}-1}\chi_{B}\circ 
T^{-s_{n,i}} - \mu(B)\big{|} d\mu
\end{align*}
and so if $\{s_{n,i}\}_{\{r_{n}\}}$ is ergodic with respect to $T$ 
then the above condition holds.  

Conversely, assume the above condition holds.  We may assume that $B$ 
is a union of 
levels.  Then for sufficiently large $n$, we may write $B = 
\bigcup_{j=0}^{\beta_{n}-1} I_{p_{n},b_{j}}$ for appropriate positive 
integers $\beta_{n}$ and $b_{j}$.  Following the techniques of Blum 
and Hanson,
{\allowdisplaybreaks
\begin{align*}
&\int\big{|}\frac{1}{r_{n}}\sum_{i=0}^{r_{n}-1}\chi_{B}\circ T^{-s_{n,i}} - 
\mu(B)\big{|}^{2} d\mu \\
&= \frac{1}{r_{n}^{2}}\sum_{i,\ell 
=0}^{r_{n}-1}\mu(T^{s_{n,\ell}}(B) \cap T^{s_{n,i}}(B)) - 
\mu(B)\mu(B) \\
&= \frac{1}{r_{n}^{2}}\sum_{i,\ell 
=0}^{r_{n}-1}\sum_{j=0}^{\beta_{n}-1}\mu(T^{s_{n,\ell}}(I_{p_{n},b_{j}}) 
\cap T^{s_{n,i}}(B)) - \mu(I_{p_{n},b_{j}})\mu(B) \\
&\leq \frac{1}{r_{n}}\sum_{\ell 
=0}^{r_{n}-1}\sum_{j=0}^{h_{p_{n}}-1}\big{|}\frac{1}{r_{n}}\sum_{i=0}^{r_{n}-1}
\mu(T^{s_{n,\ell}}(I_{p_{n},j}) \cap T^{s_{n,i}}(B)) - \mu(I_{p_{n},j}\mu(B)\big{|} \\
&\leq \frac{1}{r_{n}}\sum_{\ell 
=0}^{r_{n}-1}\Big{[}\sum_{j=0}^{h_{p_{n}}-1}\big{|}\frac{1}{r_{n}}\sum_{i=0}^{r_{n}-1}
\mu(I_{p_{n},j} \cap T^{s_{n,i}}(B)) - \mu(I_{p_{n},j}\mu(B)\big{|}
+ 2s_{n,\ell}\mu(I_{p_{n},0})\Big{]} \\
&\leq \sum_{j=0}^{h_{p_{n}}-1}\big{|}\frac{1}{r_{n}}\sum_{i=0}^{r_{n}-1}
\mu(T^{-s_{n,i}}(I_{p_{n},j}) \cap B) - \mu(I_{p_{n},j}\mu(B)\big{|} + 
\frac{2}{h_{p_{n}}}\frac{1}{r_{n}}\sum_{\ell =0}^{r_{n}-1}s_{n,\ell}.
\end{align*}}
Thus, since $\frac{1}{h_{p_{n}}r_{n}}\sum_{i=0}^{r_{n}-1}s_{n,i} \to 0$ as 
$n\to\infty$, we have that $\{s_{n,i}\}_{\{r_{n}\}}$ is ergodic with 
respect to $T$.
\end{proof}
\end{paragr}

\section{Mixing on Rank One Transformations with Restricted Growth}

\begin{paragr}\textbf{Mixing Height Sequences.}
A useful preliminary result to our main theorem are the following 
theorems equating ergodicity of the spacer sequence and mixing of the height 
sequence for any rank one transformation.
\begin{theorem}\label{T:heightmix}
Let $T$ be a rank one transformation with spacer sequence given by 
$\{s_{n,i}\}_{\{r_{n}\}}$.  Then the height sequence $\{h_{n}\}$
is mixing with respect to $T$ if and only if
$\{s_{n,i}\}_{\{r_{n}\}}$ is weak ergodic with respect to $T$.  
Equivalently, substitute the window height sequence $\{w_{n}\}$ for 
the height sequence or the representative spacer sequence $\{\widehat{s}_{n,i}\}_{\{r_{n}\}}$ 
for the spacer sequence (or both).
\end{theorem}
\begin{proof}
Let $T$, $\{s_{n,i}\}_{\{r_{n}\}}$ and $\{h_{n}\}$ be as above.  
Denote the columns defining $T$ by $C_{n}$ and the levels by $\{I_{n,j}\}$.
Let $A$ and $B$ be unions of 
levels in $C_{N}$ for some fixed $N>0$ and consider $A$ as a union of levels 
in $C_{n}$ for $n>N$: $A = \bigcup_{j=0}^{\alpha_{n}-1}I_{n,a_{j}}$ 
for some positive integers $\alpha_{n}$ and $a_{j}$ less than $h_{n}$.  
Then, 
{\allowdisplaybreaks
\begin{align*}
\mu(T^{h_{n}}&(A)\cap B) - \mu(A)\mu(B) \\
&= 
\sum_{j=0}^{\alpha_{n}-1}\sum_{i=0}^{r_{n}-1}\mu(T^{h_{n}}(I_{n,a_{j}}^{[i]})\cap B) - 
\mu(I_{n,a_{j}}^{[i]})\mu(B) \\
&= 
\sum_{j=0}^{\alpha_{n}-1}\sum_{i=0}^{r_{n}-2}\mu(T^{-s_{n,i}}(I_{n,a_{j}}^{[i+1]})\cap B) - 
\mu(I_{n,a_{j}}^{[i+1]})\mu(B) \pm \Big{(}\frac{1}{r_{n}}\Big{)} \\
&= \sum_{i=0}^{r_{n}-2}\Big{[}\sum_{\substack{j=0 \\ a_{j}\geq 
s_{n,i}}}^{\alpha_{n}-1}\mu(I_{n,a_{j}-s_{n,i}}^{[i+1]})\cap B) - 
\mu(I_{n,a_{j}-s_{n,i}}^{[i+1]})\mu(B) \\*
&\quad\quad\quad\quad + \sum_{\substack{j=0 \\ 
a_{j}<s_{n,i}}}^{\alpha_{n}-1}\mu(T^{-s_{n,i}}(I_{n,a_{j}}^{[i+1]})\cap B) - 
\mu(I_{n,a_{j}}^{[i+1]})\mu(B)\Big{]} \pm \Big{(}\frac{1}{r_{n}}\Big{)} \\
&= \sum_{i=0}^{r_{n}-2}\sum_{\substack{j=0 \\ a_{j}\geq 
s_{n,i}}}^{\alpha_{n}-1}\frac{1}{r_{n}}\big{[}\mu(I_{n,a_{j}-s_{n,i}})\cap B) - 
\mu(I_{n,a_{j}-s_{n,i}})\mu(B)\big{]} \\*
&\quad\quad\quad\quad \pm \Big{(}\sum_{i=0}^{r_{n}-2}\big{|}\mu(S_{n}^{[i]}\cap 
B) - \mu(S_{n}^{[i]})\mu(B)\big{|} + \frac{1}{r_{n}}\Big{)} \\
&= \frac{1}{r_{n}}\sum_{i=0}^{r_{n}-2}\sum_{\substack{j=0 \\ a_{j}\geq 
s_{n,i}}}^{\alpha_{n}-1}\mu(T^{-s_{n,i}}(I_{n,a_{j}})\cap B) - 
\mu(I_{n,a_{j}})\mu(B) \pm \Big{(}\mu(S_{n}) + 
\frac{1}{r_{n}}\Big{)} \\
&= 
\frac{1}{r_{n}}\sum_{i=0}^{r_{n}-1}\sum_{j=0}^{\alpha_{n}-1}\mu(T^{-s_{n,i}}(I_{n,a_{j}})\cap B) - 
\mu(I_{n,a_{j}})\mu(B) \pm \Big{(}2\mu(S_{n}) + 
\frac{2}{r_{n}}\Big{)} \\
&= \frac{1}{r_{n}}\sum_{i=0}^{r_{n}-1}\mu(T^{-s_{n,i}}(A) \cap B) - 
\mu(A)\mu(B) \pm \Big{(}2\mu(S_{n}) + \frac{2}{r_{n}}\Big{)}
\end{align*}}
which approaches zero if and only if $\{s_{n,i}\}_{\{r_{n}\}}$ is 
weak ergodic with respect to $T$ since $T$ is finite 
measure-preserving.  We may replace $\{h_{n}\}$ by $\{w_{n}\}$ since 
$\frac{w_{n}}{h_{n}} \to 1$ as $n\to\infty$ and 
$\{s_{n,i}\}_{\{r_{n}\}}$ by $\{\widehat{s}_{n,i}\}_{\{r_{n}\}}$ since 
$T$ is measure-preserving and using the $\{w_{n}\}$ substitution.
\end{proof}
\begin{theorem}\label{T:heightumix}
Let $T$ be a rank one transformation with spacer sequence given by
$\{s_{n,i}\}_{\{r_{n}\}}$.  Then the
height sequence $\{h_{n}\}$ is uniformly mixing with respect to $T$ if and only if
$\{s_{n,i}\}_{\{r_{n}\}}$ is ergodic with respect to $T$.  
Equivalently, substitute the window height sequence $\{w_{n}\}$ for 
the height sequence or the windowed spacer sequence $\{\widehat{s}_{n,i}\}_{\{r_{n}\}}$ 
for the spacer sequence (or both).
\end{theorem}
\begin{proof}
Let $T$, $\{\widehat{s}_{n,i}\}_{\{r_{n}\}}$ and $\{h_{n}\}$ be as above.  
Denote the columns defining $T$ by $C_{n}$ and the levels by 
$\{I_{n,j}\}$.  Let $B$ be a union of levels in some column $C_{N}$ 
for a fixed $N > 0$.  Then, using the same arguments as in the 
previous theorem,
{\allowdisplaybreaks
\begin{align*}
\sum_{j=0}^{h_{n}-1}&\big{|}\mu(T^{h_{n}}(I_{n,j}) \cap B) - 
\mu(I_{n,j})\mu(B)\big{|} \\
&= 
\sum_{j=0}^{h_{n}-1}\big{|}\sum_{i=0}^{r_{n}-1}\mu(T^{h_{n}}(I_{n,j}^{[i]}) 
\cap B) - \mu(I_{n,j}^{[i]})\mu(B)\big{|} \\
&= 
\sum_{j=0}^{h_{n}-1}\big{|}\frac{1}{r_{n}}\sum_{i=0}^{r_{n}-1}\mu(T^{-s_{n,i}}(I_{n,j}) \cap B)
- \mu(I_{n,j})\mu(B)\big{|} \pm o\Big{(}2\mu(S_{n}) + 
\frac{2}{r_{n}}\Big{)}
\end{align*}}
which approaches zero if and only if $\{s_{n,i}\}_{\{r_{n}\}}$ is 
ergodic with respect to $T$ by Proposition \ref{P:unifmixerg} since $T$ is 
finite measure-preserving so 
$\frac{1}{r_{n}h_{n}}\sum_{i=0}^{r_{n}-1}s_{n,i} \to 0$ as 
$n\to\infty$.  The equivalent formulations follow as above.
\end{proof}
\end{paragr}

\begin{paragr}\textbf{Mixing Sequences with Restricted Growth.}
Under the assumption that the rank one transformation in question has 
restricted growth, the transformation's mixing behavior on sequences 
is related to the ergodicity of the partial sums of the spacer 
sequence as follows; our main result follows as a consequence.
\begin{theorem}\label{T:mixseq}
Let $T$ be a restricted growth rank one transformation with spacer 
sequence $\{s_{n,i}\}_{\{r_{n}\}}$ and window height sequence $\{w_{n}\}$
and let $\{t_{n}\}$ be
a strictly increasing sequence of positive integers.  Choose 
the unique positive integers $p_{n}$ and $k_{n}$ so that $w_{p_{n}} \leq 
k_{n}h_{p_{n}} \leq t_{n} < (k_{n}+1)w_{p_{n}}$.  Denoting the 
partial sum dynamical sequences for $\{s_{n,i}\}_{\{r_{n}\}}$ by 
$\{s_{n,i}^{(k)}\}_{\{r_{n}^{(k)}\}}$, if the dynamical sequences 
$\{s_{p_{n},i}^{(k_{n})}\}_{\{r_{p_{n}}^{(k_{n})}\}}$ and 
$\{s_{p_{n},i}^{(k_{n}+1)}\}_{\{r_{p_{n}}^{(k_{n}+1)}\}}$ are both 
ergodic with respect to $T$ then $\{t_{n}\}$ is (uniform) mixing with respect 
to $T$.
\end{theorem}
\begin{proof}
Let $\{t_{m}\}$ be a strictly increasing sequence of positive integers 
and let $T$ be a restricted growth rank one 
transformation with sequence of cuts $\{r_{n}\}$ and sequence of 
spacers $\{s_{n,i}\}_{\{r_{n}\}}$.  Denote the representative spacer 
sequence for $T$ by $\{\widehat{s}_{n,i}\}_{\{r_{n}\}}$ and the family of partial sum 
dynamical sequences for the representative spacer sequence by 
$\{\widehat{s}_{n,i}^{(k)}\}_{\{r_{n}^{(k)}\}}$.  Let $\{C_{n}\}$, $\{w_{n}\}$, 
$\{I_{n,j}\}$ and $\{S_{n}\}$ be the sequences of columns, window heights, 
levels and unions of spacer levels of $T$, respectively.  To show 
uniform mixing, it suffices to show it for unions of levels; let $B$ 
be a union of levels in some column $C_{N}$ for some positive integer 
$N$.  Let $m$ be any positive integer such that $t_{m} \geq h_{N}$.  Define $p_{m}$ to be the unique 
nonnegative integer such that $w_{p_{m}} \leq t_{m} < w_{p_{m}+1}$ and let
$k_{m}$ and
$q_{m}$ be such that $t_{m}=k_{m}w_{p_{m}}+q_{m}$, where $0< 
k_{m}<r_{p_{m}}$ and $0\leq q_{m}<w_{p_{m}}$.  Since $T$ has restricted 
growth, this accounts for all sufficiently large values of $m$.  
First assume that $q_{m} < h_{p_{m}}$.
Observe that, using the techniques of the theorems on mixing height 
sequences and setting $\zeta_{n}^{(k)} = \max_{0\leq 
i<r_{n}-k}\widehat{s}_{n,i}^{(k)}$,
{\allowdisplaybreaks
\begin{align*}
&\sum_{j=0}^{h_{p_{m}}-q_{m}-1}\big{|}\mu(T^{t_{m}}(I_{p_{m},j}) \cap B) - 
\mu(I_{p_{m},j})\mu(B)\big{|} \\
&= 
\sum_{j=0}^{h_{p_{m}}-q_{m}-1}\big{|}\sum_{i=0}^{r_{p_{m}}-1}
\mu(T^{k_{m}w_{p_{m}}+q_{m}}(I_{p_{m},j}^{[i]}) \cap B)
- \mu(I_{p_{m},j}^{[i]})\mu(B)\big{|} \\
&= 
\sum_{j=0}^{h_{p_{m}}-q_{m}-1}\big{|}\sum_{i=0}^{r_{p_{m}}-k_{m}-1}
\mu(T^{k_{m}w_{p_{m}}}(I_{p_{m},j+q_{m}}^{[i]}) \cap B)
- \mu(I_{p_{m},j+q_{m}}^{[i]})\mu(B)\big{|} \\ 
&\quad  \pm \Big{(}
\sum_{j=0}^{h_{p_{m}}-q_{m}-1}\sum_{i=r_{p_{m}}-k_{m}}^{r_{p_{m}}-1}
\big{|}\mu(T^{k_{m}w_{p_{m}}}(I_{p_{m},j+q_{m}}^{[i]}) \cap B)
- \mu(I_{p_{m},j+q_{m}}^{[i]})\mu(B)\big{|}\Big{)} \\
&= 
\sum_{j=q_{m}}^{h_{p_{m}}-1}\big{|}\sum_{i=0}^{r_{p_{m}}-k_{m}-1}
\mu(T^{-\widehat{s}_{p_{m},i}^{(k_{m})}}(I_{p_{m},j}^{[i+k_{m}]}) \cap B)
- \mu(I_{p_{m},j}^{[i+k_{m}]})\mu(B)\big{|} \\ &\quad \pm \Big{(}
\sum_{j=0}^{h_{p_{m}+1}-1}\big{|}\mu(T^{h_{p_{m}+1}}(I_{p_{m}+1,j}) \cap 
B) - \mu(I_{p_{m}+1,j})\mu(B)\big{|}\Big{)} \\
&= 
\sum_{j=q_{m}}^{h_{p_{m}}-1}\big{|}\frac{1}{r_{p_{m}}}\sum_{i=0}^{r_{p_{m}}-k_{m}-1}
\mu(T^{-\widehat{s}_{p_{m},i}^{(k_{m})}}(I_{p_{m},j}) \cap B)
- \mu(I_{p_{m},j})\mu(B)\big{|} \\ &\quad \pm 
\Big{(}2\mu(I_{p_{m},0})\zeta_{n}^{(k_{m})} +
\sum_{j=0}^{h_{p_{m}+1}-1}\big{|}\mu(T^{h_{p_{m}+1}}(I_{p_{m}+1,j}) \cap 
B) - \mu(I_{p_{m}+1,j})\mu(B)\big{|}\Big{)} \\
&\leq
\sum_{j=q_{m}}^{h_{p_{m}}-1}\big{|}\frac{1}{r_{p_{m}}-k_{m}}\sum_{i=0}^{r_{p_{m}}-k_{m}-1}
\mu(T^{-\widehat{s}_{p_{m},i}^{(k_{m})}}(I_{p_{m},j}) \cap B)
- \mu(I_{p_{m},j})\mu(B)\big{|} \\ &\quad \pm 
\Big{(}\frac{2}{h_{p_{m}}}\zeta_{n}^{(k_{m})} +
\sum_{j=0}^{h_{p_{m}+1}-1}\big{|}\mu(T^{h_{p_{m}+1}}(I_{p_{m}+1,j}) \cap 
B) - \mu(I_{p_{m}+1,j})\mu(B)\big{|}\Big{)}.
\end{align*}}
Similarly, we have that
\begin{align*}
&\sum_{j=h_{p_{m}}-q_{m}}^{h_{p_{m}}-1}\big{|}\mu(T^{t_{m}}(I_{p_{m},j}) \cap B) - 
\mu(I_{p_{m},j})\mu(B)\big{|} \\
&\leq
\sum_{j=0}^{q_{m}-1}\big{|}\frac{1}{r_{p_{m}}-k_{m}-1}\sum_{i=0}^{r_{p_{m}}-k_{m}-2}\mu(T^{-\widehat{s}_{p_{m},i}^{(k_{m}+1)}}(I_{p_{m},j}) \cap B)
- \mu(I_{p_{m},j})\mu(B)\big{|} \\ &\quad \pm 
\Big{(}\frac{2}{h_{p_{m}}}\zeta_{n}^{(k_{m}+1)} +
\sum_{j=0}^{h_{p_{m}+1}-1}\big{|}\mu(T^{h_{p_{m}+1}}(I_{p_{m}+1,j}) \cap 
B) - \mu(I_{p_{m}+1,j})\mu(B)\big{|}\Big{)}.
\end{align*}
Since $T$ has restricted growth, 
$\frac{1}{h_{p_{m}}}\zeta_{n}^{(k_{m})} \to 0$ as 
$m\to\infty$ and similarly for $k_{m}+1$.
If $\{s_{n,i}\}_{\{r_{n}\}}$ is uniformly ergodic with respect to $T$ 
then it is ergodic with respect to $T$ so $\{h_{n}\}$ is a uniformly 
mixing sequence with respect to $T$ by Theorem \ref{T:heightumix}.  
By Proposition \ref{P:unifmixerg}, the remaining above quantities then approach zero as $m\to\infty$.  Thus, $T$ is 
uniformly mixing so $T$ is mixing by Proposition \ref{P:unifmix}. 

When $q_{m} \geq h_{p_{m}}$, set $q_{m}^{\prime} = q_{m} - h_{p_{m}}$ 
and note that $\frac{q_{m}^{\prime}}{h_{p_{m}}}\to 0$ since $T$ is 
finite measure-preserving.  Then, as above,
{\allowdisplaybreaks
\begin{align*}
&\sum_{j=0}^{h_{p_{m}}-q_{m}^{\prime}-1}\big{|}\mu(T^{m}(I_{p_{m},j}) \cap B) - 
\mu(I_{p_{m},j})\mu(B)\big{|} \\
&\leq
\sum_{j=q_{m}^{\prime}}^{h_{p_{m}}-1}\big{|}\frac{1}{r_{p_{m}}-k_{m}-1}\sum_{i=0}^{r_{p_{m}}-k_{m}-2}
\mu(T^{-\widehat{s}_{p_{m},i}^{(k_{m}+1)}+\bar{s}_{n}}(I_{p_{m},j}) \cap B)
- \mu(I_{p_{m},j})\mu(B)\big{|} \\ &\quad \pm 
\Big{(}\frac{2}{h_{p_{m}}}\zeta_{n}^{(k_{m}+1)} +
\sum_{j=0}^{h_{p_{m}+1}-1}\big{|}\mu(T^{h_{p_{m}+1}}(I_{p_{m}+1,j}) \cap 
B) - \mu(I_{p_{m}+1,j})\mu(B)\big{|}\Big{)}.
\end{align*}}
Then, also as above, using Theorem \ref{T:heightumix} and Proposition 
\ref{P:unifmixerg}, the quantities above approach zero when 
$\{s_{n,i}\}_{\{r_{n}\}}$ is uniformly ergodic with respect to $T$. 
\end{proof}
\end{paragr}

\begin{paragr}\textbf{Mixing Transformations with Restricted Growth.}
Our main result follows from the previous theorem relating mixing on 
sequences to ergodic averages of the partial sums of the spacer 
sequence.
\begin{theorem}\label{T:main}
Let $T$ be a restricted growth rank one transformation.  Then $T$ is 
a mixing transformation if and only if the spacer sequence for $T$ is 
uniformly ergodic.
\end{theorem}
\begin{proof}
Let $T$ be as above.  If the spacer sequence for $T$ is uniformly 
ergodic, then by Theorem~\ref{T:mixseq}, every strictly increasing 
sequence of positive integers is mixing with respect to $T$ so $T$ is 
mixing.
Conversely, 
if $\{s_{n,i}\}_{\{r_{n}\}}$ is not uniformly ergodic with respect 
to $T$ then along some sequence of positive integers $\{k_{n}\}$, the 
dynamical sequence $\{s_{n,i}^{(k_{n})}\}_{\{r_{n}^{(k_{n})}\}}$ is not 
ergodic with respect to $T$.  We may assume that $\frac{k_{n}}{r_{n}}$ 
is bounded away from zero and one since the (standard) ergodicity of 
the spacer sequence and the fact that $T$ is rank one would then imply 
the result.  Then there exists a $\delta > 0$ and a 
union of levels $B$ in some column $C_{N}$ for some positive integer 
$N$ such that
\[
\int\big{|}\frac{1}{r_{n}^{(k_{n})}}\sum_{i=0}^{r_{n}^{(k_{n})}-1}\chi_{B}\circ 
T^{-s_{n,i}^{(k_{n})}} - \mu(B)\big{|}^{2} \geq \delta.
\]
For any positive integer $n > N$, write $B$ as a union of levels in 
$C_{n}$; $B = \bigcup_{j=0}^{\beta_{n}-1}I_{n,b_{j}}$.  Using the 
techniques above,
{\allowdisplaybreaks
\begin{align*}
&\frac{1}{(r_{n}-k_{n})^{2}}\sum_{i,j=0}^{r_{n}-k_{n}-1}
\mu(T^{s_{n,i}^{(k_{n})}-s_{n,j}^{(k_{n})}}(B) \cap B) - 
\mu(B)\mu(B) \\
&= \sum_{z=0}^{\beta_{n}-1}\frac{1}{(r_{n}-k_{n})^{2}}\sum_{i,j=0}^{r_{n}-k_{n}-1}
\mu(T^{s_{n,i}^{(k_{n})}-s_{n,j}^{(k_{n})}}(I_{n,b_{z}}) \cap B) - 
\mu(I_{n,b_{z}})\mu(B) \\
&= 
\frac{r_{n}}{r_{n}-k_{n}}\sum_{z=0}^{\beta_{n}-1}\frac{1}{r_{n}-k_{n}}\sum_{i,j=0}^{r_{n}-k_{n}-1}
\mu(T^{s_{n,i}^{(k_{n})}-s_{n,j}^{(k_{n})}}(I_{n,b_{z}}^{[j+k_{n}]}) 
\cap B) - \mu(I_{n,b_{z}}^{[j+k_{n}]})\mu(B) \\
&= 
\frac{r_{n}}{r_{n}-k_{n}}\sum_{z=0}^{\beta_{n}-1}\frac{1}{r_{n}-k_{n}}\sum_{i,j=0}^{r_{n}-k_{n}-1}
\mu(T^{s_{n,i}^{(k_{n})}+k_{n}h_{n}}(I_{n,b_{z}}^{[j]}) 
\cap B) - \mu(I_{n,b_{z}}^{[j]})\mu(B) \\ &\quad \pm 
\Big{(}\frac{1}{h_{p_{m}}}(\max_{i}\widehat{s}_{p_{m},i}^{(k_{n})})\Big{)} \\
&= \frac{r_{n}}{r_{n}-k_{n}}\sum_{z=0}^{\beta_{n}-1}\frac{1}{r_{n}-k_{n}}\sum_{i=0}^{r_{n}-k_{n}-1}
\mu(T^{s_{n,i}^{(k_{n})}+k_{n}h_{n}}(I_{n,b_{z}}) 
\cap B) - \mu(I_{n,b_{z}})\mu(B) \\ &\quad \pm 
\Big{(}\frac{1}{h_{p_{m}}}(\max_{i}\widehat{s}_{p_{m},i}^{(k_{n})}) + 
\sum_{j=0}^{h_{p_{m}+1}-1}\big{|}\mu(T^{h_{p_{m}+1}}(I_{p_{m}+1,j}) \cap 
B) - \mu(I_{p_{m}+1,j})\mu(B)\big{|}\Big{)} \\
&= \frac{r_{n}}{r_{n}-k_{n}}\sum_{z=0}^{\beta_{n}-1}\frac{1}{r_{n}-k_{n}}\sum_{i=0}^{r_{n}-k_{n}-1}
\mu(T^{s_{n,i}^{(k_{n})}+k_{n}h_{n}}(B) \cap B) - \mu(B)\mu(B) \\
&\quad \pm 
\Big{(}\frac{1}{h_{p_{m}}}(\max_{i}\widehat{s}_{p_{m},i}^{(k_{n})}) + 
\sum_{j=0}^{h_{p_{m}+1}-1}\big{|}\mu(T^{h_{p_{m}+1}}(I_{p_{m}+1,j}) \cap 
B) - \mu(I_{p_{m}+1,j})\mu(B)\big{|}\Big{)}
\end{align*}}
and so since $\{h_{n}\}$ is a uniformly mixing sequence and $T$ has 
restricted growth, the dynamical sequence $\{s_{n,i}^{(k_{n})} + 
k_{n}h_{n}\}_{\{r_{n}-k_{n}\}}$ is not weak ergodic with respect to 
$T$.  Since $\{s_{n,i}^{(k_{n})} + 
k_{n}h_{n}\}_{\{r_{n}-k_{n}\}}$ is weak monotone, 
by the generalized weak Blum-Hanson theorem 
(Theorem~\ref{T:genweakBH}), this means that $T$ is not mixing.
\end{proof}
\end{paragr}

\section{Staircase Transformations}\label{S:appl}

\begin{paragr}\textbf{Construction of Staircase Transformations.}
The class of staircase transformations has appeared in the 
literature recently providing examples of rank one mixing 
transformations.  We 
include this section to demonstrate the application of our theorems to 
explicit rank one 
constructions, yielding an alternate proof of the result shown in \cite{tA98}.
Formally, a rank one transformations $T$ with cut sequence $\{r_{n}\}$ and 
spacer sequence $\{s_{n,i}\}_{\{r_{n}\}}$ is a 
\textbf{\emph{staircase transformation}} when the spacers are given 
by $s_{n,i} = i$ (a ``staircase'' pattern) for all $n$ and $0 \leq i 
< r_{n}$.
By Proposition \ref{P:rnbounded}, if the sequence $\{r_{n}\}$ has a 
finite limit point, then $T$ cannot be mixing; we assume
that staircase transformations have an unbounded cut sequence.
Note that 
restricted growth on staircase transformations is equivalent to 
$\frac{r_{n}^{2}}{h_{n}}\to 0$ as $n\to\infty$ where $\{r_{n}\}$ is 
the sequence of cuts and $\{h_{n}\}$ is the sequence of heights for 
the staircase.  This is precisely the condition given by Adams in 
\cite{tA98} for staircase transformations to be mixing.
\end{paragr}

\begin{paragr}\textbf{Weak Mixing on Staircase Transformations.}
The first step in showing mixing is showing weak mixing, accomplished 
by the following application of our preliminary result on mixing 
height sequences.
\begin{theorem}\label{T:stairmixheight}
Let $T$ be a staircase transformation.  Then $T$ has a mixing height 
sequence; hence, $T$ is a weak mixing transformation.
\end{theorem}
\begin{proof}
Let $T$ be a staircase transformation with spacer sequence 
$\{s_{n,i}\}_{\{r_{n}\}}$.  Then the spacer sequence is given by 
$s_{n,i} = i$ is ergodic with respect to $T$ since $T$ is ergodic,
so $T$ has a mixing height sequence and is thus weak 
mixing so $T$ has a weak mixing spacer sequence.
\end{proof}
\end{paragr}

\begin{paragr}\textbf{Uniform Weak Mixing on Staircase Spacer 
Sequences.}
We begin with two lemmas, due to Adams \cite{tA98}.  The first 
follows directly from the measure-preserving property; the proof is left to the reader.
\begin{lemma}\label{L:pLR}(\cite{tA98})\textbf{\emph{Block Lemma}}
Let $T$ be a measure-preserving transformation and $B$ a measurable 
set.  Then for any positive integers $R$, $L$ and $p$,
\[
\int \big{|}\frac{1}{R}\sum_{i=0}^{R-1}\chi_{B}\circ T^{-i} - \mu(B)\big{|} d\mu 
\leq \int \big{|}\frac{1}{L}\sum_{i=0}^{L-1}\chi_{B}\circ T^{-ip} - \mu(B)\big{|} 
d\mu + \frac{pL}{R}.
\]
\end{lemma}
\begin{lemma}\label{L:BHtech}(\cite{tA98})
Let $T$ be an ergodic transformation and $\{v_{n}\}$ a sequence of 
positive integers such that for any fixed (nonzero) integer $i$, the 
sequence of integers $\{iv_{n}\}$ is mixing with respect to $T$.  
Then for any $\epsilon > 0$ there exist arbitrarily 
large positive integers $L$ such that for sufficiently large $n$,
\[
\int \big{|}\frac{1}{L}\sum_{i=0}^{L-1}\chi_{B\times B}\circ (T\times 
T)^{-iv_{n}} - 
\mu\times\mu(B\times B)\big{|} d\mu\times\mu < \epsilon.
\]
\end{lemma}
\begin{proof}
Fix $\epsilon > 0$.  For each fixed nonzero integer $i$, choose a
positive integer
$N_{i}$ such that for all $n \geq N_{i}$, $\big{|}\mu(T^{iv_{n}}(B) \cap 
B) - \mu(B)\mu(B)\big{|} < \epsilon$.  For any positive integer $L$, 
set $N = \sup_{-L\leq i\leq L}N_{i}$.  Note that for any $0 < j,\ell \leq 
L$, $-L \leq j - \ell \leq L$. Then, for all $n \geq N$, first 
applying the H\"{o}lder Inequality,
{\allowdisplaybreaks
\begin{align*}
\Big{[}\int &\big{|}\frac{1}{L}\sum_{i=0}^{L-1}\chi_{B\times B}\circ (T\times T)^{-iv_{n}} - 
\mu\times\mu(B\times B)\big{|} d\mu\times\mu\Big{]}^{2} \\
&\leq
\int \big{|}\frac{1}{L}\sum_{j=1}^{L}\chi_{B\times B}\circ (T\times 
T)^{-jv_{n}} - 
\mu\times\mu(B\times B)\big{|}^{2} d\mu\times\mu \\ &= 
\frac{1}{L^{2}}\sum_{j,\ell=1}^{L}\mu\times\mu((T\times T)^{(j-\ell 
)v_{n}}(B\times B) \cap B\times B) - 
\mu\times\mu(B\times B)^{2} \\ &< 
\frac{1}{L^{2}}\sum_{j,\ell =1}^{L}\epsilon^{2} + 
\mu(B)\mu(B)\epsilon = \epsilon(\epsilon +\mu(B)\mu(B)).
\end{align*}}
\end{proof}
\begin{proposition}\label{P:stairuwm}
Let $T$ be a staircase transformation.  Then $T$ has a uniformly weak 
mixing spacer sequence.
\end{proposition}
\begin{proof}
Let $T$ be a staircase transformation with spacer sequence 
$\{s_{n,i}\}_{\{r_{n}\}}$ and let the partial sum dynamical sequences 
for $\{s_{n,i}\}_{\{r_{n}\}}$ be denoted
$\{s_{n,i}^{(k)}\}_{\{r_{n}^{(k)}\}}$.  Note that $T$ is weak mixing by 
Theorem~\ref{T:stairmixheight}.  The partial sum 
dynamical sequences for the spacer sequence are given by $s_{n,i}^{(k)} = 
ik + \frac{1}{2}k(k+1)$.  Hence, for each fixed positive integer $k$, 
the dynamical sequence takes on the maximum value 
$\widetilde{s}_{n}^{(k)} = k(r_{n}-k-1) + \frac{1}{2}k(k-1)$ so has 
density $D_{n}^{(k)} = \frac{r_{n}-k_{n}}{\widetilde{s}_{n}^{(k)}} \to
\frac{1}{k}$ as $n\to\infty$.  Theorem~\ref{T:weakmixdens} then 
implies that $\{s_{n,i}^{(k)}\}_{\{r_{n}^{(k)}\}}$ is weak mixing with 
respect to $T$ for each fixed $k$.

Let $\{k_{n}\}$ be any 
sequence of positive integers such that $k_{n} < r_{n}$ for all $n$ 
and $0 < \liminf_{n\to\infty}\frac{k_{n}}{r_{n}} = \gamma$.  
Note that $\frac{h_{n}^{2}}{h_{n+1}}\to\infty$ as 
$n\to\infty$, since, on a finite space, $\frac{h_{n}}{r_{n}} \to \infty$ and 
$\frac{r_{n}h_{n}}{h_{n+1}} \to 1$ as $n \to \infty$, so
$\frac{h_{n}^{2}}{h_{n+1}} = 
\frac{r_{n}h_{n}}{h_{n+1}}\frac{h_{n}}{r_{n}} \to \infty$.
Define the 
sequence $\{p_{n}\}$ so that $h_{p_{n}-1}\leq k_{n} < h_{p_{n}}$ for 
each $n$.  Define $\{u_{n}\}$ by $u_{n} = 
\inf \{u\in\mathbb{Z} : uk_{n}\geq h_{p_{n}}\}$.  Then, $h_{p_{n}}\leq 
u_{n}k_{n} < 2h_{p_{n}}$ so 
\begin{align*}
\frac{u_{n}}{r_{n}-k_{n}} &= 
\frac{u_{n}k_{n}}{k_{n}(r_{n}-k_{n})} < 
\frac{2h_{p_{n}}}{k_{n}k_{n}}\frac{k_{n}}{r_{n}-k_{n}} \\ &\leq 
\frac{2h_{p_{n}}}{h_{p_{n}-1}^{2}}\frac{k_{n}}{r_{n}}\frac{r_{n}}{r_{n}-k_{n}} \leq 
2(1-\gamma)\frac{h_{p_{n}}}{h_{p_{n}-1}^{2}}
\end{align*}
which approaches zero as $n\to\infty$ since 
$\frac{h_{n-1}^{2}}{h_{n}} \to \infty$ as $n\to\infty$.
For any fixed nonzero integer $i$, 
consider the sequence of integers 
$\{ik_{n}u_{n}\}$ is mixing with respect to $T$.  Clearly, 
$ih_{p_{n}} \leq ik_{n}u_{n} < 2ih_{p_{n}}$.  Set $j_{n}$ to be the 
positive integer such that $ik_{n}u_{n} = j_{n}h_{p_{n}}$.  Then 
$i\leq j_{n} < 2i$ for all $n$ so the dynamical sequences given by 
$\{s_{n,z}^{(j_{n})}\}_{\{r_{n}^{(j_{n})}\}}$ and 
$\{s_{n,z}^{(j_{n}+1)}\}_{\{r_{n}^{(j_{n}+1)}\}}$ are weak mixing 
with respect to $T$ since $i$ is fixed.  Theorem~\ref{T:mixseq} then implies that the 
sequence $\{ik_{n}u_{n}\}$ is mixing with respect to $T$.

Let $\{A_{n}\}$ be any sequence of measurable sets and $B$ be any 
measurable set.  Fix $\epsilon > 0$.  Lemma~\ref{L:BHtech} yields 
positive integers $L$ and $N$ such that for all integers $n \geq N$,
\[
\int \big{|}\frac{1}{L}\sum_{i=0}^{L-1}\chi_{B\times B}\circ (T\times T)^{-i\rho_{n}} - 
\mu\times\mu(B\times B)\big{|} d\mu\times\mu < \epsilon.
\]
Since $\frac{u_{n}}{r_{n}-k_{n}} \to 0$ as $n\to\infty$, there exists 
an integer $\widehat{N} \geq N$ such that for all integers 
$n\geq\widehat{N}$, $L\frac{u_{n}}{r_{n}-k_{n}} < \epsilon$.  Then, 
using the block lemma (Lemma~\ref{L:pLR}),
{\allowdisplaybreaks
\begin{align*}
\frac{1}{r_{n}^{(k_{n})}}&\sum_{i=0}^{r_{n}^{(k_{n})}-1}\big{|}
\mu(T^{s_{n,i}^{(k_{n})}}(A_{n}) \cap B) - 
\mu(A_{n})\mu(B)\big{|}^{2} \\
&= \int_{A_{n}\times A_{n}}\frac{1}{r_{n}^{(k_{n})}}\sum_{i=0}^{r_{n}^{(k_{n})}-1}
\chi_{B\times B}\circ (T\times T)^{s_{n,i}^{(k_{n})}} - 
\mu\times\mu(B\times B) d\mu\times\mu \\
&\leq \int\big{|}\frac{1}{r_{n}^{(k_{n})}}\sum_{i=0}^{r_{n}^{(k_{n})}-1}
\chi_{B\times B}\circ (T\times T)^{s_{n,i}^{(k_{n})}} - 
\mu\times\mu(B\times B)\big{|} d\mu\times\mu \\
&\leq \int \big{|}\frac{1}{L}\sum_{i=0}^{L-1}\chi_{B\times B}\circ (T\times T)^{-i\rho_{n}} - 
\mu\times\mu(B\times B)\big{|} d\mu\times\mu + 
L\frac{u_{n}}{r_{n}-k_{n}} \\ &< 2\epsilon.
\end{align*}}
\end{proof}
\end{paragr}

\begin{paragr}\textbf{Mixing on Restricted Growth Staircase Transformations.}
Using our above results and our main theorem, we prove the following 
result originally in \cite{tA98}.
\begin{theorem}[Adams]
Let $T$ be a staircase transformation that has restricted growth.  
Then $T$ is a mixing transformation.
\end{theorem}
\begin{proof}
Proposition~\ref{P:stairuwm} and Theorem~\ref{T:main}.
\end{proof}
\end{paragr}

\section{Ornstein's ``Random Spacers'' Method}\label{S:randomspacers}

\begin{paragr}\textbf{Construction with ``Random Spacers''.}
We conclude with a discussion of the mixing rank one transformations due 
to Ornstein \cite{dO72} using a ``random spacers'' method for cutting 
and stacking.  The reader is referred to \cite{mN98} for a detailed 
account of this method.  The transformations are defined by choosing a 
set of values $\{x_{n,i}\}_{i=0}^{r_{n}-1}$ using the uniform 
distribution on the set of integers between $-\frac{1}{2}\widetilde{s}_{n}$ 
and $\frac{1}{2}\widetilde{s}_{n}$ where $\{\widetilde{s}_{n}\}$ is a given 
sequence of positive integers with no finite limit points (Ornstein's original constructions 
used $\widetilde{s}_{n} = h_{n-1}$) and letting the spacer sequence for the 
transformation $T$ be given by $s_{n,i} = \widetilde{s}_{n} + x_{n,i+1} - 
x_{n,i}$ where $x_{n,r_{n}} = x_{n,0}$.  The sequence of cuts $\{r_{n}\}$ 
is a sequence of positive integers specified later to show mixing properties.  
The window height sequence $\{w_{n}\}$ is then given, letting 
$\{h_{n}\}$ denote the height sequence, by 
$w_{n} = h_{n} + \big{\lfloor}\widetilde{s}_{n} + \frac{x_{n,r_{n}} - 
x_{n,0}}{r_{n}}\big{\rfloor} = h_{n} + \widetilde{s}_{n}$.  For any 
positive integers $n$, $i$ and $k$ such that $0\leq i < i + k < r_{n}$, the partial sum of 
the representative spacer sequence
$\big{|}\widehat{s}_{n,i}^{(k)}\big{|} = \big{|}x_{n,i+k} - x_{n,i}\big{|} 
\leq \widetilde{s}_{n}$ so the 
transformation $T$ has restricted growth since 
$\frac{\widetilde{s}_{n}}{h_{n}} \to 0$ as $n\to\infty$ is a necessary 
condition for $T$ to be finite measure-preserving.
\end{paragr}

\begin{paragr}\textbf{Probabilistic Lemma.}
To show mixing, we will need the following Lemma used by Ornstein \cite{dO72}; the 
proof may be found in \cite{mN98}.
\begin{lemma}\label{L:prob}
Let $H$ be a positive integer and $X = \{i\in\mathbb{Z}:\big{|} 
i\big{|}\leq\frac{H}{2}\}$.  For any positive integer $m$, let 
$\Omega_{m} = X^{m}$ and let $P_{m}$ be the uniform distribution on 
$\Omega_{m}$.  Let $\omega = (\omega_{1},\ldots ,\omega_{m}) \in 
\Omega_{m}$ and let $x_{i}$, $0 < i \leq m$, denote the coordinates 
of the random variable on $\Omega_{m}$.  For each integer $\ell$, set 
$C_{k,\ell} = \#\{0<i\leq m : x_{i+k}(\omega) - x_{i}(\omega) = 
\ell\}$.  Then given $\alpha > 1$, $\epsilon > 0$ and a positive 
integer $N$, there exists an integer $m \geq N$ such that
\[
P_{m}\big{(}\bigcap_{k=1}^{(1-\epsilon)m}\bigcap_{\ell\in\mathbb{Z}}\{\omega : 
C_{k,\ell} \leq \frac{\alpha}{H}(m-k)\}\big{)} > 1 - \epsilon.
\]
\end{lemma}
\end{paragr}

\begin{paragr}\textbf{Weak Mixing using ``Random Spacers''.}
We first show that almost surely such a transformation is weak 
mixing when the spacer sequence has positive upper density.  Note that El Houcein has shown that 
almost surely such transformations are totally ergodic without 
requiring our condition on the sequence of cuts and sequence of ranges; the reader is 
referred to \cite{mN98}.
Our proof is accomplished using the techniques in \cite{BFMS01}; the reader is referred to 
that work for details on double ergodicity.
\begin{theorem}
Let $T$ be a rank one transformation constructed using ``random 
spacers'' as above with sequence of cuts $\{r_{n}\}$ and sequence of ranges for 
the spacer sequence $\{\widetilde{s}_{n}\}$ such that 
$\limsup_{n\to\infty}\frac{r_{n}}{\widetilde{s}_{n}} > 0$---the 
spacer sequence for $T$ has positive (upper) density.  
Then almost surely $T$ is a weak mixing transformation.
\end{theorem}
\begin{proof}
Let $T$, $\{r_{n}\}$, and $\{\widetilde{s}_{n}\}$ be as above and let $\epsilon > 
0$.  Let $A$ and $B$ be measurable sets with $\mu(A) > 0$ and $\mu(B) > 
0$.  Then there exist levels $I$ and $J$ in some 
column $C_{N}$ for some positive integer $N$ such that $I$ and $J$ 
are $(1-\epsilon)$-full of $A$ and $B$, respectively.  Let $\ell$ be 
the distance between $I$ and $J$ in $C_{N}$ ($\ell$ positive when $I$ 
is above $J$).  Then $\mu(I \cap A) + \mu(J \cap B) > 
2(1-\epsilon)\mu(I)$.  Choose $\gamma > 0$ such that $\mu(I \cap A) + 
\mu(J \cap B) - 2(1-\epsilon)\mu(I) \geq 8\epsilon\gamma\mu(I) > 0$.  
For any integer $n$ such that $n \geq N$ and $\frac{2}{r_{n}} \leq 
\gamma$, write $I$ and $J$ as unions of levels in $C_{n+1}$; denote 
the sublevels of $I$ in $C_{n+1}$ by $I_{t}$ for $0 \leq t < R_{n}$ 
where $R_{n} = \prod_{z=N}^{n}r_{z}$.  Order the $I_{t}$ so that each 
block of sublevels $I_{mr_{n}+k}$ for $0\leq k <r_{n}$ and $m$ fixed forms a 
whole level in $C_{n}$ such that $I_{mr_{n}+k}$ is in the $k$th 
subcolumn of $C_{n}$.  Similarly, we have $J_{t}$ for $0 \leq t < 
R_{n}$ ordered such that $I_{t}$ is $\ell$ above $J_{t}$ where $\ell$ 
is the distance between $I$ and $J$ as above.  Then
\[
\sum_{k=0}^{R_{n}-1}\mu(I_{t}\cap A) + \mu(I_{t+1}\cap A) + 
2\mu(J_{t+2}\cap B) - 4(1-\epsilon )\mu(I_{t}) \geq 
16\epsilon\gamma\mu(I).
\]
Let $x_{n}$ denote the number of values of $t$ for $0 \leq t < R_{n}$ 
with $t\mod r_{n} \ne r_{n}-1$ and $t\mod r_{n} \ne r_{n}-2$
such that $\mu(I_{t}\cap A) + \mu(I_{t+1}\cap A) + 2\mu(J_{t+2}\cap 
B) - 4(1-\epsilon)\mu(I_{t}) > 0$.  Then
\begin{align*}
\sum_{t=0}^{R_{n}-1}&\mu(I_{t}\cap A) + \mu(I_{t+1}\cap A) + 
2\mu(J_{t+2}\cap B) - 4(1-\epsilon )\mu(I_{t}) \\
&\leq (R_{n} - x_{n}) 0 + \big{(}\frac{2R_{n}}{r_{n}} + 
x_{n}\big{)}\big{(}4\mu(I_{t}) - 4(1-\epsilon )\mu(I_{t})\big{)} \\
&= \big{(}\frac{2R_{n}}{r_{n}} + 
x_{n}\big{)}4\epsilon\frac{1}{R_{n}}\mu(I) = 
\frac{8\epsilon\mu(I)}{r_{n}} + \frac{4\epsilon x_{n}\mu(I)}{R_{n}}.
\end{align*}
Thus $16\epsilon\gamma\mu(I) \leq \frac{8\epsilon\mu(I)}{r_{n}} + 
\frac{4\epsilon x_{n}\mu(I)}{R_{n}}$ so $\frac{4\epsilon 
x_{n}\mu(I)}{R_{n}} \geq 8\epsilon\mu(I)\big{(}2\gamma - 
\frac{2}{r_{n}}\big{)} \geq 8\epsilon\mu(I)(2\gamma - \gamma)$; 
hence $x_{n} \geq 2\gamma R_{n}$.  Consider the possible values of 
$t\mod r_{n}$ for the at least $\gamma R_{n}$ values of $t$ such that 
$\mu(I_{t}\cap A) + \mu(I_{t+1}\cap A) + 2\mu(J_{t+2}\cap 
B) - 4(1-\epsilon)\mu(I_{t}) > 0$.  At most $R_{n-1}$ choices for $t$ have the 
same value (mod $r_{n}$) so there must exist at least $2\gamma r_{n}$ 
distinct values of $k$ with $0 \leq k < r_{n}$ such that 
$\mu(I_{t}\cap A) + \mu(I_{t+1}\cap A) + 2\mu(J_{t+2}\cap 
B) - 4(1-\epsilon)\mu(I_{t}) > 0$ where $t = m r_{n} + k$ for 
some $m$.  Observe 
that since $P_{n}$ is the uniform distribution on the integers 
between $-\frac{\widetilde{s}_{n}}{2}$ and $\frac{\widetilde{s}_{n}}{2}$, using the notation from Lemma~\ref{L:prob},
\[
P_{n}\big{(}\{\omega : s_{n,k} - s_{n,k+1} \ne \ell \quad\text{for all 
$\lfloor\gamma r_{n}\rfloor$ values $k$}\}\big{)} = \big{(}1 - 
\frac{1}{\widetilde{s}_{n}}\big{)}^{\lfloor\gamma r_{n}\rfloor}.
\]
Thus
{\allowdisplaybreaks
\begin{align*}
P\big{(}&\bigcup_{n=N}^{\infty}\{\omega : s_{n,k} - s_{n,k+1} = \ell 
\quad\text{for some $k$ of the $\lfloor\gamma r_{n}\rfloor$ 
values}\}\big{)} \\
&= 1 - P\big{(}\bigcap_{n=N}^{\infty}\{\omega : s_{n,k} - s_{n,k+1} 
\ne \ell \quad\text{for all $\lfloor\gamma r_{n}\rfloor$ 
values $k$}\}\big{)} \\
&= 1 - \prod_{n=N}^{\infty}P_{n}\big{(}\{\omega : s_{n,k} - s_{n,k+1} 
\ne \ell \quad\text{for all
$\lfloor\gamma r_{n}\rfloor$ values $k$}\}\big{)} \\
&= 1 - \prod_{n=N}^{\infty}\big{(}1 - 
\frac{1}{\widetilde{s}_{n}}\big{)}^{\lfloor\gamma r_{n}\rfloor}.
\end{align*}}
Using the approximation $\log(1-x) \approx -x$ for small $x$, we have 
that
{\allowdisplaybreaks
\begin{align*}
\log\prod_{n=N}^{\infty}\big{(}1-\frac{1}{\widetilde{s}_{n}}\big{)}^{\lfloor\gamma r_{n}\rfloor} &=
\sum_{n=N}^{\infty}\lfloor\gamma 
r_{n}\rfloor\log\big{(}1-\frac{1}{\widetilde{s}_{n}}\big{)} \\
&\approx 
\sum_{n=N}^{\infty}\gamma r_{n}\frac{-1}{\widetilde{s}_{n}} = 
-\gamma\sum_{n=N}^{\infty}\frac{r_{n}}{\widetilde{s}_{n}}.
\end{align*}}
By our requirement that $\frac{r_{n}}{\widetilde{s}_{n}}$ be bounded away from 
zero along some subsequence, this implies that 
\[
\prod_{n=N}^{\infty}\big{(}1-\frac{1}{\widetilde{s}_{n}}\big{)}^{\lfloor\gamma r_{n}\rfloor} = 0.
\]
Hence, we have that
\[
P\big{(}\bigcup_{n=N}^{\infty}\{\omega : s_{n,k} - s_{n,k+1} = \ell 
\quad\text{for some $k$ of the $\lfloor\gamma r_{n}\rfloor$ 
values}\}\big{)} = 1.
\]
Therefore, almost surely there exists $I_{t}$, $I_{t+1}$ and $J_{t+2}$ 
such that $\mu(I_{t} \cap A) + \mu(I_{t+1} \cap A) + 2\mu(J_{t+2} \cap B) > 
4(1-\epsilon)\mu(I_{t})$ and $s_{n,k} - s_{n,k+1} = \ell$ for 
some $n$ and $k < r_{n} - 2$ where $t = m r_{n} + k$ for some $m$.  In this case, we have that
\begin{align*}
T^{h_{n}+s_{n,k}}(I_{t}) &= I_{t+1}; \\
T^{h_{n}+s_{n,k}}(I_{t+1}) &= 
T^{s_{n,k}-s_{n,k+1}}(I_{t+2}) = T^{\ell}(I_{t+2}) = 
J_{t+2}.
\end{align*}
Since $\mu(I_{t} \cap A) > (1-4\epsilon)\mu(I_{t})$, 
$\mu(I_{t+1} \cap A) > (1-4\epsilon)\mu(I_{t+1})$ and 
$\mu(J_{t+2} \cap B) > (1-4\epsilon)\mu(J_{t+2})$, 
we then have that
{\allowdisplaybreaks
\begin{align*}
\mu(T^{h_{n}+s_{n,k}}(A) \cap A) &\geq 
(1-4\epsilon)\mu(T^{h_{n}+s_{n,k}}(I_{t}) \cap A) \\ &= 
(1-4\epsilon)\mu(I_{t+1} \cap A) \geq 
(1-4\epsilon)^{2}\mu(I_{t+1}) > 0
\end{align*}}
and, similarly,
{\allowdisplaybreaks
\begin{align*}
\mu(T^{h_{n}+s_{n,k}}(A) \cap B) &\geq 
(1-4\epsilon)\mu(T^{h_{n}+s_{n,k}}(I_{t+1}) \cap B) \\ &= 
(1-4\epsilon)\mu(J_{t+2} \cap B) \geq 
(1-4\epsilon)^{2}\mu(J_{t+2}) > 0.
\end{align*}}
Thus, $T$ is doubly ergodic which is equivalent to weak mixing.
\end{proof}
\end{paragr}

\begin{paragr}\textbf{Mixing using ``Random Spacers''.}
Assume that the transformation $T$ has been partially constructed up 
to the column $C_{n-1}$ using the ``random spacers'' method.  Apply Lemma~\ref{L:prob} 
with $H = \widetilde{s}_{n}$, 
$N = r_{n-1}$, a fixed $\alpha > 1$ and an $\epsilon_{n} > 0$ such 
that $\epsilon_{n} \to 0$ as $n\to\infty$ to obtain $r_{n} = m$.  Set 
$C_{k,\ell}^{n} = \#\{0 < i \leq r_{n} - k : x_{i+k}(\omega) - 
x_{i}(\omega) = \ell\}$ and
\[
L_{n} = \bigcap_{k=1}^{(1-\epsilon_{n})r_{n}}\bigcap_{\ell = 
-\frac{1}{2}\widetilde{s}_{n}}^{\frac{1}{2}\widetilde{s}_{n}}\{\omega : C_{k,\ell}^{n} \leq 
\frac{\alpha}{\widetilde{s}_{n}}(r_{n} - k)\}
\]
so that $P_{r_{n}}(L_{n}) > 1 - \epsilon_{n}$.  Using $\omega \in 
L_{n}$, the representative spacers for $T$ are given by 
$\widehat{s}_{n,i} = x_{n,i+1}(\omega) - x_{n,i}(\omega) + \widetilde{s}_{n} - 
\bar{s}_{n} = x_{n,i+1}(\omega) - x_{n,i}(\omega)$ so the 
partial sums of the representative spacer sequence are given by 
$\widehat{s}_{n,i}^{(k)} = x_{n,i+k}(\omega) - x_{n,i}(\omega)$ and 
$r_{n}^{(k)} = r_{n} - k$.  Then, for any $0 < k < r_{n} - 1$ and any 
fixed measurable set $B$ and any sequence of measurable sets $\{A_{n}\}$,
{\allowdisplaybreaks
\begin{align*}
\frac{1}{r_{n}^{(k)}}\sum_{i=0}^{r_{n}^{(k)}-1}&\big{|}
\mu(T^{\widehat{s}_{n,i}^{(k)}}(A_{n}) \cap B) - \mu(A_{n})\mu(B)\big{|} \\
&= \frac{1}{r_{n}-k}\sum_{\ell = 
-\frac{1}{2}\widetilde{s}_{n}}^{\frac{1}{2}\widetilde{s}_{n}}C_{k,\ell}^{n}\big{|}\mu(T^{\ell}(A_{n}) \cap B) - 
\mu(A_{n})\mu(B)\big{|} \\
&\leq \frac{1}{r_{n}-k}\sum_{\ell = 
-\frac{1}{2}\widetilde{s}_{n}}^{\frac{1}{2}\widetilde{s}_{n}}\frac{\alpha}{\widetilde{s}_{n}}(r_{n} - 
k)\big{|}\mu(T^{\ell}(A_{n}) \cap B) - \mu(A_{n})\mu(B)\big{|} \\
&= \frac{\alpha}{\widetilde{s}_{n}}\sum_{\ell = 
-\frac{1}{2}\widetilde{s}_{n}}^{\frac{1}{2}\widetilde{s}_{n}}\big{|}\mu(T^{\ell}(A_{n}) \cap B) 
- \mu(A_{n})\mu(B)\big{|}.
\end{align*}}
Since $\widetilde{s}_{n} \to \infty$ as $n \to \infty$ and $T$ is weak mixing, 
for any sequence of positive integers $\{k_{n}\}$ such that $k_{n} < 
r_{n} - 1$ for all $n$, the dynamical 
sequence $\{\widehat{s}_{n,i}^{(k_{n})}\}_{\{r_{n}^{(k_{n})}\}}$ is weak 
mixing with respect to $T$ so is ergodic with respect to $T$ by 
Proposition~\ref{P:wmimperg}.  Theorem \ref{T:main} then yields the 
following theorem originally in \cite{dO72}.
\begin{theorem}[Ornstein]
Let $T$ be a rank one transformation constructed using ``random 
spacers'' with sequence of cuts increasing sufficiently fast as above.
Then almost surely $T$ is a mixing transformation.
\end{theorem}
\end{paragr}


\begin{thebibliography}{99}

\bibitem[Ad98]{tA98}T. Adams, {\it Smorodinsky's Conjecture on Rank One Mixing},
  Proc. Amer. Math. Soc. 126 (1998), 739-744. 
  
\bibitem[AF92]{AF92}T. Adams and N. Friedman, \emph{Staircase Mixing}, unpublished (1992).
  
\bibitem[AS99]{AS99}T. Adams and C. Silva, {\it $\mathbb{Z}^{d}$-Staircase Actions},
  Ergodic Th. and Dyn. Sys. 19 (1999), 837-850.    

\bibitem[BH60]{BH60}J. Blum and D. Hanson, {\it On the Mean Ergodic Theorem for Subsequences},
  Bull. Amer. Math. Soc. 55 (1960), 308-311.
  
\bibitem[Bo93]{jB93}J. Bourgain, {\it On the Spectral Type of Ornstein's Class One Transformations}, 
  Israel J. Math. 84 (1993), no. 1-2, 53--63. 
  
\bibitem[BFMS01]{BFMS01}A. Bowles, L. Fidkowski, A. Marinello and C. 
  Silva, {\it Double Ergodicity of Nonsingular Transformations and 
  Infinite Measure Preserving Staircase Transformations}, Illinois J. 
  Math (to appear).

\bibitem[Fe97]{sF97}S. Ferenczi, {\it Systems of Finite Rank},
  Colloq. Math. 73 (1997), 35-65.
  
\bibitem[Fr70]{nF70}N. Friedman, {\bf Introduction to Ergodic Theory},
  Van Nostrand (1970).
  
\bibitem[Fr83]{nF83}N. Friedman, {\it Mixing on Sequences}, Can. J. 
  Math, 35 (1983), 339-352.

\bibitem[Fu81]{hF81}H. Furstenberg, {\bf Recurrence in Ergodic Theory and Combinatorial 
  Number Theory}, Princeton University Press, Princeton (1981).
  
\bibitem[Ho99]{eA99}El A. El Houcein,
  {\it La singularitŽ mutuelle presque s\^{u}r du spectre des transformations d'Ornstein}
  Israel J. Math. 112 (1999), 135--155.
  
\bibitem[Ka84]{sK84}S.  Kalikow, {\it Twofold Mixing Implies 
  Threefold Mixing for Rank One Transformations},
  Ergodic Th. and Dyn. Sys. 4 (1984), no. 2, 237--259. 

\bibitem[Ki88]{jK88}J. King, {\it Joining-rank and the Structure of Finite 
  Rank Mixing Transformations}, J. Analyse Math. 51 (1988), 182--227.

\bibitem[Kl96]{iK96}I. Klemes, {\it The Spectral Type of the Staircase Transformation},
  Tohoku Math. J. (2) 48 (1996), no. 2, 247--258. 
  
\bibitem[Na98]{mN98}M. Nadkarni, {\bf Spectral Theory of Dynamical 
  Systems}, Birkh\"{a}user (1998).
 
\bibitem[Or72]{dO72}D. Ornstein, {\it On the Root Problem in Ergodic Theory}, 
  Proc. of the Sixth Berkeley Symposium on Mathematical Statistics and 
  Probability, Univ. of California Press (1972), 347-356.

\bibitem[Ry93]{vR93}V. Ryzhikov, {\it Mixing, Rank and Minimal 
  Self-Joinings of Actions with an Invariant Measure}, Russian Acad. 
  Sci. Sb. Math, Vol. 75, No. 2 (1993), 405-426.

\end{thebibliography}
\end{document}